\documentclass[12pt]{amsart}
\usepackage[pagebackref=true,colorlinks,linkcolor=blue,citecolor=magenta]{hyperref}
\usepackage{amsmath,amssymb,latexsym} 
\usepackage{graphicx}
\usepackage{fullpage}
\usepackage{enumitem}
\usepackage{comment}
\usepackage{etoolbox}
\usepackage[square, numbers,sort&compress]{natbib}
\usepackage{hyperref}
\usepackage{cleveref}

\usepackage{algorithmicx,algcompatible,algorithm}

\usepackage{multirow}

\usepackage{tikzSettings}

\theoremstyle{plain}
\newtheorem{theo}{Theorem}
\newtheorem{prop}[theo]{Proposition}

\newtheorem{lem}[theo]{Lemma}

\theoremstyle{remark}
\newtheorem{rem}{Remark}

\crefname{ru}{rule}{rules}
\Crefname{ru}{Rule}{Rules}
\crefname{tip}{principle}{principles}
\Crefname{tip}{Principle}{Principles}

\newenvironment{outdent}
{\begin{list}{}{\leftmargin-2cm\rightmargin\leftmargin}\centering\item\relax}
{\end{list}\ignorespacesafterend}

\newenvironment{changemargin}[2]{\begin{list}{}{%
\setlength{\topsep}{0pt}%
\setlength{\leftmargin}{0pt}%
\setlength{\rightmargin}{0pt}%
\setlength{\listparindent}{\parindent}%
\setlength{\itemindent}{\parindent}%
\setlength{\parsep}{0pt plus 1pt}%
\addtolength{\leftmargin}{#1}%
\addtolength{\rightmargin}{#2}%
}\item }{\end{list}}

\AtBeginEnvironment{proof-claim}{\vspace{-0.5cm}\footnotesize\begin{changemargin}{0.5cm}{0.5cm}}
\AtEndEnvironment{proof-claim}{\end{changemargin}}

\newcommand{\ds}{\displaystyle}
\newcommand{\ve}{\boldsymbol}

\usepackage{xspace}

\newcommand{\MRCSP}{\textsc{Monoid Resource Constrained Shortest Path Problem}\xspace}

\newcommand{\pairing}{p}
\newcommand{\pairingSet}{\mathcal{P}}

\newcommand{\leg}{\ell}
\newcommand{\legSet}{\mathcal{L}}

\newcommand{\maintenance}{\Delta_{\text{maint}}}						

\newcommand{\lpProp}{\alpha}

\newcommand{\ldProp}{\beta}

\newcommand{\scConSet}{\mathsf{S}}
\newcommand{\scSet}{S}

\newcommand{\re}{q}
\newcommand{\rleq}{\preceq}

\newcommand{\rplus}{\oplus}
\newcommand{\bigrplus}{\bigoplus}
\newcommand{\rset}{M}
\newcommand{\rcost}{c}
\newcommand{\rmeas}{\rho}
\newcommand{\meet}{\wedge} 

\newcommand{\join}{\vee} 

\newcommand{\calR}{\mathcal{R}}

\def\UB{\operatorname{UB}}

\def\L{\mathcal{L}}
\def\R{\mathbb{R}}
\def\Z{\mathbb{Z}}

\def\C{\mathcal{C}}

\def\bb{\boldsymbol{b}}
\def\cc{\boldsymbol{c}}

\def\xx{\boldsymbol{x}}
\def\yy{\boldsymbol{y}}
\def\zz{\boldsymbol{z}}
\newcommand{\bfx}{\xx}
\newcommand{\bfy}{\yy}
\newcommand{\bfz}{\zz}
\newcommand{\ovbfx}{\overline{\xx}}
\newcommand{\ovbfy}{\overline{\yy}}
\newcommand{\ovbfz}{\overline{\zz}}
\newcommand{\ovx}{\overline{x}}
\newcommand{\ovy}{\overline{y}}
\newcommand{\ovz}{\overline{z}}

\newcommand{\ind}{\mathbf{1}}

\usepackage[normalem]{ulem}

\usepackage{todonotes}

\newcommand{\na}{n^{\mathrm{a}}} 

\title{Aircraft routing and crew pairing: \\ updated algorithms at Air France}

\author{Axel Parmentier}
\address{A. Parmentier,
Universit\'e Paris Est, CERMICS, 77455 Marne-la-Vall\'ee CEDEX, France}
\email{axel.parmentier@enpc.fr}
\author{Fr\'ed\'eric Meunier}
\address{F. Meunier, Universit\'e Paris Est, CERMICS, 77455 Marne-la-Vall\'ee CEDEX, France}
\email{frederic.meunier@enpc.fr}


\date{\today}

\begin{document}

\maketitle

\begin{abstract}
Aircraft routing and crew pairing problems aim at building the sequences of flight legs operated respectively by airplanes and by crews of an airline. Given their impact on airlines operating costs, both have been extensively studied for decades. Our goal is to provide reliable and easy to maintain frameworks for both problems at Air France. We propose simple approaches to deal with Air France current setting. For routing, we introduce an exact compact IP formulation that can be solved to optimality by current MIP solvers in at most a few minutes even on Air France largest instances. Regarding crew pairing, we provide a  methodology to model the column generation pricing subproblem within a new resource constrained shortest path framework recently introduced by the first author. This new framework, which can be used as a black-box, leverages on bounds to discard partial solutions and speed-up the resolution. The resulting approach enables to solve to optimality Air France largest instances. Recent literature has focused on integrating aircraft routing and crew pairing problems. As a side result, we are able to solve to near optimality large industrial instances of the integrated problem by combining the aforementioned algorithms within a simple cut generating method. 
\end{abstract}



\section{Introduction}

\subsection{Context}
Interactions between Operations Research and Air Transport Industry have been successful for at least five decades \cite{barnhart2003applications,etschmaier1974operations}. These interactions have taken various forms: yield management, airplane timetabling, ground operations scheduling, air traffic management, etc. Key applications are notably the construction of sequences of flight legs operated by airplanes and crews. As airplane sequences of legs are \emph{routes} and crew sequences \emph{pairings}, this construction is called \emph{aircraft routing} for airplanes, and \emph{crew pairing} for crews. 


The present paper focuses on these two applications and is the fruit of a research partnership with Air France, the main French airline. 
We aim at providing a reliable and easy to maintain framework that can cope with the specific and challenging industrial context of the company. Air France working rules are 
more complex 
than the IATA standards: Collective agreements reach hundreds of pages, and two of the most cited references on crew pairing \citep{desaulniers1997crew,lavoie1988new} develop ad-hoc approaches to build pairings satisfying the company's rules. While the aircraft routing remains easy, this turns the exact resolution of Air France crew pairing into a challenge.


As crews need time to cross airports if they change airplane, the two problems are linked and the sequential resolution currently in use in the industry is suboptimal.
 Solving the integrated problem has been identified by academics as a difficult problem. Air France also requested an easy to maintain solution scheme for the integrated problem.


\subsection{Literature review}\label{subsec:review} All the versions of the problems considered in this paper are NP-hard. We focus here on mathematical programming approaches.

Aircraft routing is considered at Air France as a pure feasibility problem, which contrasts with the recent literature which considers optimization versions. Authors either maximize profit when the fleet is heterogeneous \citep{desaulniers1997daily,barnhart1998flight}, or minimize delay propagation along sequences of flights \citep{lan2006planning}. 
A recent paper introduces tools to deal with richer maintenance constraints \citep{safaei2018aircraft}.
State-of-the-art solution approaches rely on column generation \citep{barnhart1998flight,desaulniers1997daily,lan2006planning,froyland2013recoverable}, where columns are sequences of flight legs between airports where maintenance checks can be performed. 
They can solve to optimality large instances of the optimization versions in a few hours. (The solution proposed in \cite{lan2006planning} is actually for the so-called {\em tail assignment problem}, where the airplanes are distinguishable, but it can be adapted to aircraft routing.)
 Alternative approaches include heuristics \citep{feo1989flight} and Lagrangian relaxations~\citep{clarke1997aircraft}. 

\citet{clarke1997aircraft} propose a MIP with few variables. However, their MIP has two exponential-size families of constraints: The first one is formed by classical subtour elimination constraints; the second one is formed by ``minimal violation path'' constraints that enforce the maintenance requirements. While the first family can be discarded for the aircraft routing problem we consider, the second one must be kept, and thus cut generation cannot be avoided to solve their MIP. Practically efficient compact integer programming formulations, that is, formulations which require neither column generation, nor cut generation, have recently been proposed by Cacchiani and Salazar-Gonz\'alez \citep{cacchiani2016optimal} and \citet{khaled2018compact}. Compact formulations have the advantage to be more handy and can often be directly implemented in standard MIP solvers. \citet{cacchiani2016optimal} consider the integrated problem, and propose a compact integer programming formulation for aircraft routing. Since they assume that the airplanes spend alternatively one night in a base and one night outside, their approach does not generalize to the Air France case, where one maintenance must be performed at least every four days, while the maintenance day is not fixed. 
By the way, due to the different maintenance requirements, the aircraft routing problem considered by \citet{cacchiani2016optimal} has a polynomial status \citep{gopalan1998aircraft}, while Air France problem is NP-hard~\citep{talluri1998four}.
\citet{khaled2018compact} have also recently proposed a compact MIP approach to tail assignment, which can be adapted to aircraft routing. They are able to solve to optimality instances with up to 1,178 legs and 30 airplanes in 3 hours. Their approach could in principle be used to address the version of the aircraft routing problem met at Air France. We discuss later in the paper the advantage of the approach we propose for aircraft routing with respect to theirs.

Crew variable wages and hotel rooms are among airlines first sources of variable costs. As both depend on the sequences of flight legs crews operate, crew pairing is an intensively studied optimization problem; see \citet{gopalakrishnan2005airline} for an extensive review. 
Regulatory agencies and collective bargaining agreements list numerous working rules that make the crew pairing problem highly non-linear and hence difficult. 
There is a long tradition of MIP approaches to crew pairing \citep{gopalakrishnan2005airline,baker1981computational}.
Since the seminal work of \citet{minoux1984column}, state-of-the-art approaches solve the crew pairing by column generation \cite{lavoie1988new,anbil1992global,hoffman1993solving,chu1997solving,desaulniers1997crew,klabjan2002airline,klabjan1999airline,barnhart1998approximate,subramanian2008effective,zeren2016novel}. They consider set partitioning formulations where columns are possible pairings. These methods hide the non-linearity in the pricing subproblem, which can be efficiently solved using resource constrained shortest path approaches \citep{irnich2005shortest}.  As a large part of the working rules apply to \emph{duties}, i.e., 
subsequences of a pairing formed by the flight legs operated on a same day, the subproblem is often split into two parts \citep{desaulniers1997crew,zeren2016novel}. 
The first one builds the set of all non-dominated duties. The second one builds the pairings by solving a path problem in the graph whose vertices are the non-dominated duties, and whose arcs are the pair of duties that can be chained.
However, as the number of non-dominated duties is huge, solving the pricing subproblem becomes costly on large instances.
When working rules are simple, one can also use compact integer programming approaches \citep{beasley1996tree} where variables indicate if a given connection is used, and set partitioning formulations, where columns are the duties \cite{vance1997airline}. However, this is generally not the case, and such models are generally turned into initialization heuristics \cite{barnhart1998approximate}.

During the last decade, much attention has been devoted to the integration of aircraft routing and crew pairing. Moving from a sequential to an integrated approach enables to reduce the cost by $5\%$ on average according to \citet{cordeau2001benders}, and $1.6\%$ according to \citet{papadakos2009integrated}.   
Solution methods are column generation based heuristics \cite{mercier2005computational,mercier2007integrated,cohn2003improving,papadakos2009integrated,weide2010iterative,salazar2014approaches,cacchiani2016optimal,shao2015novel}. 
The heuristics of \citet{weide2010iterative,salazar2014approaches,cacchiani2016optimal} share many similarities with the ones we propose in this paper. Dunbar et al.~consider robust~\citep{dunbar2012robust} and stochastic~\citep{dunbar2014integrated} versions of the problem. 
To the best of our knowledge, 
the largest instances considered in the literature have 750 legs \citep{weide2010iterative}.

\subsection{Contribution and methods}

The present paper is the result of a project initiated by Air France to design efficient and easy to maintain solution schemes for aircraft routing, crew pairing, and the integrated problem.

The first author has recently proposed an abstract framework~\cite{parmentier2015algorithms} for computing resource constrained shortest paths. The main contribution of the present paper is the proof that this framework can be used on a concrete problem and considerably improves the size of the instances that can be solved at optimality. Indeed, we apply this framework to the pricing subproblem of a standard column generation approach for the crew pairing problem and solve to optimality from a few minutes to a few hours instances with up to $1,000$ flight legs, which outperforms previous performances on that problem. One key element in the performances of this framework is the use of sets of bounds to discard paths, instead of single bounds: This is useful in a context where any two resources are not necessarily comparable (in Section~\ref{sec:pricing_subproblem}, further details will be given). Even if this idea of sets of bounds is present in the aforementioned paper of the first author, the present paper is the first proof that such a technique is very efficient in practice. We finally emphasize that the framework for shortest path computation does not explain how to model concrete problems like the one met for crew pairing.
The modeling we propose is thus also a contribution on its own: More than $70$ rules have to be satisfied, and most of them are non-linear. 
Finally, Desrosiers and L\"ubbecke \citep[p.~16]{desaulniers2006column} underline that, in a column generation context, ``accelerating the pricing algorithm itself usually leads most significant speeds-up''. As all the approaches to the integrated problem use a column generation approach for the crew pairing, we believe that these approaches can be significantly accelerated by using our improved crew pairing pricing subproblem algorithm.


Our second contribution is a simple and compact integer formulation for aircraft routing. 
Such formulations are desirable in an industrial context since, as mentioned in the literature review, they do not require tricky development and can often be directly implemented in off-the-shelf solvers. In addition, our formulation is very efficient: It enables to solve all Air France industrial instances in at most a few minutes.
As mentioned in the literature review,  other compact formulations have recently been proposed \cite{cacchiani2016optimal,khaled2018compact}. 
Even though the formulation of Khaled et al.~\citep{khaled2018compact} could in theory be adapted to the Air France specific problem, 
such an adaptation is not straightforward due to the fact that routes are cyclic in the Air France aircraft routing problem.
Furthermore, on the Air France problem, our formulation admits a stronger linear relaxation than theirs (we discuss it in Appendix~\ref{sec:aircraft_routing_mip}), which gives a clear competitive advantage to our approach.
We do not claim that our formulation outperforms their one on other versions of the problem such as the ones they consider.


Finally, we design a simple  cut generating method for solving the integrated problem, which relies on our contributions for aircraft routing and crew pairing.
Like the one of \citet{weide2010iterative}, our method consists in solving alternatively crew pairing and aircraft routing problems. However, they do not consider the same problem: Their version includes a notion of robustness with respect to delay.
Experiments show that the method is able to solve to near optimality instances with up to $1,766$ flight legs, which again outperforms previous results on that problem.  Due to the specificity of the Air France problem, with no aircraft routing costs, our method cannot handle all the problems considered in the literature. However, it is the only one that proves optimality gaps smaller than 0.01\% on instances with more than 600 legs.

We emphasize that for the three problems, our solution is easy to use and to maintain by the company. The algorithm for computing the shortest paths is already implemented and can be used as a black-box. The only non-trivial task is the modeling of the rules in the framework, but once a few techniques have been understood (like the ones we use later in the paper, in Section~\ref{sub:modeling_the_pricing_subproblem}), even this step is straightforward.

\subsection{Organization of the paper}

Each of the three problems considered in the paper is addressed in a separate section: The aircraft routing problem is studied in Section~\ref{sec:ar}, the crew pairing problem in Section~\ref{sec:crewPairing}, and the integrated problem in Section~\ref{sec:integrated_problem}. Each of these sections gets exactly the same structure. It starts with a subsection describing the problem. A second subsection is then devoted to a modeling of the problem (e.g., the compact integer formulation for the aircraft routing problem). It ends with a subsection explaining the proposed method to solve the problem (e.g., column generation for the crew pairing problem and cut generation for the integrated problem).

Experiments showing the efficiency of the methods proposed in each of these three sections are provided and discussed in Section~\ref{sec:experimentalResults}. 

Section~\ref{sec:pricing_subproblem} is a section making a focus on an algorithmic subroutine required by our method for the crew pairing problem. This section is much more technical than the others and can be safely skipped at first reading (and the same holds for Section~\ref{sub:focus_on_the_pricing_subproblem_choice_of_} that deals with specific experiments regarding this subroutine). This subroutine is the algorithm solving the pricing subproblem of the column generation. It relies on the shortest path framework of the first author and on bound sets, both described in that section.

The paper ends with a short conclusion (Section~\ref{sec:conclusion}). All proofs are postponed to Appendix~\ref{sec:proofs}.

%

\section{Compact integer program for aircraft routing}\label{sec:ar}

\subsection{Problem formulation}\label{subsec:ar_formulation}

\begin{figure}
\begin{center}
\begin{tikzpicture}
\def\l{1.5};
\def\h{1.3};
\node[vert, fill=blue, fill opacity=0.2 ,text opacity=1] (jb) at (0*\l,2*\h) {$\ell_l$};
\node[vert] (a)  at (1*\l,2*\h) {$\ell_a$};
\node[vert] (b)  at (2*\l,2*\h) {$\ell_b$};
\node[vert] (c)  at (3*\l,2*\h) {$\ell_c$};
\node[vert] (d)  at (4*\l,2*\h) {$\ell_d$};
\node[vert, fill = blue, fill opacity=0.2, text opacity=1] (e)  at (6*\l,1.2*\h) {$\ell_e$};
\node[vert, fill = blue, fill opacity=0.2, text opacity=1] (eb) at (0*\l,1.2*\h) {$\ell_f$};
\node[vert] (f)  at (1*\l,1.2*\h) {$\ell_g$};
\node[vert] (g)  at (2*\l,1.2*\h) {$\ell_h$};
\node[vert] (h)  at (3*\l,1.2*\h) {$\ell_i$};
\node[vert] (i)  at (5*\l,1.2*\h) {$\ell_j$};
\node[vert, fill=blue, fill opacity=0.2 ,text opacity=1] (j)  at (6*\l,2*\h) {$\ell_k$};
\node[vert] (l)  at (1*\l,0*\h) {$\ell_m$};
\node[vert] (m)  at (2*\l,0*\h) {$\ell_n$};
\node[vert] (n)  at (3*\l,0*\h) {$\ell_o$};
\node[vert] (o)  at (4*\l,0*\h) {$\ell_p$};
\node[vert, fill=purple, fill opacity=0.2, text opacity=1] (k)  at (6*\l,0*\h) {$\ell_q$};
\node[vert, fill=purple, fill opacity=0.2, text opacity=1] (kb)  at (0*\l,0*\h) {$\ell_r$};

\draw[->,>=latex, blue] (a) -- (b);
\draw[->,>=latex, blue] (b) -- (c);
\draw[->,>=latex, blue] (c) -- (d);
\draw[->,>=latex, blue] (d) -- (e);
\draw[->,>=latex, blue] (f) -- (g);
\draw[->,>=latex, blue] (g) -- (h);
\draw[->,>=latex, blue] (h) -- (i);
\draw[->,>=latex, blue] (i) -- (j);
\draw[->,>=latex, blue] (jb) -- (a);
\draw[->,>=latex, blue] (eb) -- (f);
\draw[->,>=latex,blue] (j) -- ++(0,.7) -- ++(-6*\l,0) -- (jb);
\draw[->,>=latex,blue] (e) -- ++(0,-0.7) -- ++(-6*\l,0) -- (eb);

\draw[dashed, ->,>=latex, red] (kb) -- (l);
\draw[dashed, ->,>=latex, red] (l) -- (m);
\draw[dashed, ->,>=latex, red] (m) -- (n);
\draw[dashed, ->,>=latex, red] (n) -- (o);
\draw[dashed, ->,>=latex, red] (o) -- (k);
\draw[dashed, ->,>=latex, red] (k) -- ++(0,-0.7) -- ++(-6*\l,0) -- (kb);

\draw (0.5*\l,3*\h) -- (0.5*\l,-0.8*\h);
\node  at (0*\l,2.8*\h) {Day 7};
\node  at (1*\l,2.8*\h) {Day 1};
\draw (2.5*\l,3*\h) -- (2.5*\l,-0.8*\h);
\node  at (3*\l,2.8*\h) {Day 2};

\draw (5.5*\l,3*\h) -- (5.5*\l,-0.8*\h);
\node  at (6*\l,2.8*\h) {Day 7};

\draw[dotted] (4.5*\l,3*\h) -- (4.5*\l,-0.8*\h);

\end{tikzpicture}
\end{center}
\caption{Two routes. A two-week route in plain line. A single-week route in dashed line}
\label{fig:routeExamples}
\end{figure}


Building the sequences of flight legs required for aircraft routing corresponds to solving the following problem.

The input is formed by a set of airports, a collection $\L$ of {\em flight legs}, and a number $n^{\mathrm{a}}$ of airplanes. Some airports are {\em bases} in which maintenance checks can be performed. A flight leg is characterized by departure and arrival airports, as well as departure and arrival times (it is of course assumed that departure time is smaller than arrival time for any flight leg). We consider the flight legs on a weekly horizon: the departure and arrival times are given for a typical week.

The purpose of aircraft routing is to determine routes for airplanes so that each flight leg is operated by an airplane each week without using more than  $n^{\mathrm{a}}$ airplanes. In addition, there are maintenance operations that have to be regularly performed: Each airplane has to spend a night in a base at least every $\maintenance$ days, where $\maintenance$ is a given parameter, which is equal to $4$ at Air France.

Formally, an \emph{airplane connection} is a pair $(\ell,\ell')$ of flight legs which satisfies
\begin{itemize}
\item the arrival airport of $\ell$ is the departure airport of $\ell'$
\item the duration between the departure time of $\ell'$ and the arrival time of $\ell$ is bounded from below by a fixed quantity (which can depend on the airport, the time, and the fleet).
\end{itemize}
We underline that there are connections $(\ell,\ell')$ with $\ell$ at the end of the week and $\ell'$ at the beginning of the (next) week.
A {\em route} is a cyclic sequence of distinct flight legs $\ell_1,\ldots,\ell_k$ such that any two consecutive flight legs $(\ell,\ell') = (\ell_{i-1},\ell_i)$ or $(\ell_k,\ell_1)$ is an \emph{airplane connection}.
Routes can last several weeks, but each week, the sequences of flight legs operated by airplanes are the same. In other words, when we consider all airplanes as indiscernible, the solution must have a week periodicity. Figure~\ref{fig:routeExamples}~illustrates two routes, which last respectively one and two weeks. 
As each flight leg has to be operated each week, routes lasting $p$ weeks require $p$ airplanes. A route satisfies the \emph{maintenance requirement} if an airplane following this route in a cyclic way (repeating the solution when it reaches the end of the cycle) spends a night in a base at least every $\maintenance$ days

The task consists in partitioning $\L$ into routes satisfying the maintenance requirement such that the number of airplanes needed to operate these routes is less than or equal to $n^{\mathrm{a}}$.



\subsection{Modeling as an integer program}\label{subsec:ip-ar}

\begin{figure}
\begin{center}
\begin{tikzpicture}

\def\l{1.5}
\def\h{0.5}

\draw (-0.5*\l,5*\h) -- (9.5*\l,5*\h) node (e) {};

\node[above left = 0.3 and -0.5 of e, text width=2cm, align=center] {Legs not ending in a base};
\node[below left = 0.3 and -0.5 of e, text width=2cm, align=center] {Legs ending in a base};

\draw (0.5*\l,-1.5*\h) -- (0.5*\l,13*\h); 
\draw (1.5*\l,-1.5*\h) -- (1.5*\l,13*\h); 
\draw (3.5*\l,-1.5*\h) -- (3.5*\l,13*\h); 
\draw (4.5*\l,-1.5*\h) -- (4.5*\l,13*\h); 
\draw (5.5*\l,-1.5*\h) -- (5.5*\l,13*\h); 
\draw (6.5*\l,-1.5*\h) -- (6.5*\l,13*\h); 
\draw (7.5*\l,-1.5*\h) -- (7.5*\l,13*\h); 

\node at (0*\l,12.5*\h) {Day 7};
\node at (1*\l,12.5*\h) {Day 1};
\node at (2.5*\l,12.5*\h) {Day 2};
\node at (4*\l,12.5*\h) {Day 3};
\node at (5*\l,12.5*\h) {Day 4};
\node at (6*\l,12.5*\h) {Day 5};
\node at (7*\l,12.5*\h) {Day 6};
\node at (8*\l,12.5*\h) {Day 7};

\node (as) at (0,5.7*\h) {};
\node[above = 0*\h of as] (a0) {$\circ$};
\node[above = 1*\h of as] (a1) {$\circ$};
\node[above = 2*\h of as] (a2) {$\circ$};
\node[above = 3*\h of as] (a3) {$\circ$};
\node[above = 4*\h of as] (an) {$\ell_a$};
\node[rounded corners,draw,fit=(as.north)(an)] {};

\node (cs) at (1*\l,5.7*\h) {};
\node[above = 0*\h of cs] (c0) {$\circ$};
\node[above = 1*\h of cs] (c1) {$\circ$};
\node[above = 2*\h of cs] (c2) {$\circ$};
\node[above = 3*\h of cs] (c3) {$\circ$};
\node[above = 4*\h of cs] (cn) {$\ell_c$};
\node[rounded corners,draw,fit=(cs.north)(cn)] {};

\node (es) at (2*\l,5.7*\h) {};
\node[above = 0*\h of es] (e0) {$\circ$};
\node[above = 1*\h of es] (e1) {$\circ$};
\node[above = 2*\h of es] (e2) {$\circ$};
\node[above = 3*\h of es] (e3) {$\circ$};
\node[above = 4*\h of es] (en) {$\ell_e$};
\node[rounded corners,draw,fit=(es.north)(en)] {};

\node (gs) at (3*\l,5.7*\h) {};
\node[above = 0*\h of gs] (g0) {$\circ$};
\node[above = 1*\h of gs] (g1) {$\circ$};
\node[above = 2*\h of gs] (g2) {$\circ$};
\node[above = 3*\h of gs] (g3) {$\circ$};
\node[above = 4*\h of gs] (gn) {$\ell_g$};
\node[rounded corners,draw,fit=(gs.north)(gn)] {};

\node (is) at (4*\l,5.7*\h) {};
\node[above = 0*\h of is] (i0) {$\circ$};
\node[above = 1*\h of is] (i1) {$\circ$};
\node[above = 2*\h of is] (i2) {$\circ$};
\node[above = 3*\h of is] (i3) {$\circ$};
\node[above = 4*\h of is] (in) {$\ell_i$};
\node[rounded corners,draw,fit=(is.north)(in)] {};

\node (ks) at (5*\l,5.7*\h) {};
\node[above = 0*\h of ks] (k0) {$\circ$};
\node[above = 1*\h of ks] (k1) {$\circ$};
\node[above = 2*\h of ks] (k2) {$\circ$};
\node[above = 3*\h of ks] (k3) {$\circ$};
\node[above = 4*\h of ks] (kn) {$\ell_k$};
\node[rounded corners,draw,fit=(ks.north)(kn)] {};

\node (ms) at (6*\l,5.7*\h) {};
\node[above = 0*\h of ms] (m0) {$\circ$};
\node[above = 1*\h of ms] (m1) {$\circ$};
\node[above = 2*\h of ms] (m2) {$\circ$};
\node[above = 3*\h of ms] (m3) {$\circ$};
\node[above = 4*\h of ms] (mn) {$\ell_m$};
\node[rounded corners,draw,fit=(ms.north)(mn)] {};

\node (os) at (7*\l,5.7*\h) {};
\node[above = 0*\h of os] (o0) {$\circ$};
\node[above = 1*\h of os] (o1) {$\circ$};
\node[above = 2*\h of os] (o2) {$\circ$};
\node[above = 3*\h of os] (o3) {$\circ$};
\node[above = 4*\h of os] (on) {$\ell_o$};
\node[rounded corners,draw,fit=(os.north)(on)] {};

\node (qs) at (8*\l,5.7*\h) {};
\node[above = 0*\h of qs] (q0) {$\circ$};
\node[above = 1*\h of qs] (q1) {$\circ$};
\node[above = 2*\h of qs] (q2) {$\circ$};
\node[above = 3*\h of qs] (q3) {$\circ$};
\node[above = 4*\h of qs] (qn) {$\ell_q$};
\node[rounded corners,draw,fit=(qs.north)(qn)] {};


\node (bs) at (0,-1.5*\h) {};
\node[above = 0*\h of bs] (b0) {$\circ$};
\node[above = 1*\h of bs] (b1) {$\circ$};
\node[above = 2*\h of bs] (b2) {$\circ$};
\node[above = 3*\h of bs] (b3) {$\circ$};
\node[above = 4*\h of bs] (bn) {$\ell_b$};
\node[rounded corners,draw,fit=(bs.north)(bn)] {};

\node (ds) at (1*\l,-1.5*\h) {};
\node[above = 0*\h of ds] (d0) {$\circ$};
\node[above = 1*\h of ds] (d1) {$\circ$};
\node[above = 2*\h of ds] (d2) {$\circ$};
\node[above = 3*\h of ds] (d3) {$\circ$};
\node[above = 4*\h of ds] (dn) {$\ell_d$};
\node[rounded corners,draw,fit=(ds.north)(dn)] {};

\node (fs) at (2*\l,-1.5*\h) {};
\node[above = 0*\h of fs] (f0) {$\circ$};
\node[above = 1*\h of fs] (f1) {$\circ$};
\node[above = 2*\h of fs] (f2) {$\circ$};
\node[above = 3*\h of fs] (f3) {$\circ$};
\node[above = 4*\h of fs] (fn) {$\ell_f$};
\node[rounded corners,draw,fit=(fs.north)(fn)] {};

\node (hs) at (3*\l,-1.5*\h) {};
\node[above = 0*\h of hs] (h0) {$\circ$};
\node[above = 1*\h of hs] (h1) {$\circ$};
\node[above = 2*\h of hs] (h2) {$\circ$};
\node[above = 3*\h of hs] (h3) {$\circ$};
\node[above = 4*\h of hs] (hn) {$\ell_h$};
\node[rounded corners,draw,fit=(hs.north)(hn)] {};

\node (js) at (4*\l,-1.5*\h) {};
\node[above = 0*\h of js] (j0) {$\circ$};
\node[above = 1*\h of js] (j1) {$\circ$};
\node[above = 2*\h of js] (j2) {$\circ$};
\node[above = 3*\h of js] (j3) {$\circ$};
\node[above = 4*\h of js] (jn) {$\ell_j$};
\node[rounded corners,draw,fit=(js.north)(jn)] {};

\node (ls) at (5*\l,-1.5*\h) {};
\node[above = 0*\h of ls] (l0) {$\circ$};
\node[above = 1*\h of ls] (l1) {$\circ$};
\node[above = 2*\h of ls] (l2) {$\circ$};
\node[above = 3*\h of ls] (l3) {$\circ$};
\node[above = 4*\h of ls] (ln) {$\ell_l$};
\node[rounded corners,draw,fit=(ls.north)(ln)] {};

\node (ns) at (6*\l,-1.5*\h) {};
\node[above = 0*\h of ns] (n0) {$\circ$};
\node[above = 1*\h of ns] (n1) {$\circ$};
\node[above = 2*\h of ns] (n2) {$\circ$};
\node[above = 3*\h of ns] (n3) {$\circ$};
\node[above = 4*\h of ns] (nn) {$\ell_n$};
\node[rounded corners,draw,fit=(ns.north)(nn)] {};

\node (ps) at (7*\l,-1.5*\h) {};
\node[above = 0*\h of ps] (p0) {$\circ$};
\node[above = 1*\h of ps] (p1) {$\circ$};
\node[above = 2*\h of ps] (p2) {$\circ$};
\node[above = 3*\h of ps] (p3) {$\circ$};
\node[above = 4*\h of ps] (pn) {$\ell_p$};
\node[rounded corners,draw,fit=(ps.north)(pn)] {};

\node (rs) at (8*\l,-1.5*\h) {};
\node[above = 0*\h of rs] (r0) {$\circ$};
\node[above = 1*\h of rs] (r1) {$\circ$};
\node[above = 2*\h of rs] (r2) {$\circ$};
\node[above = 3*\h of rs] (r3) {$\circ$};
\node[above = 4*\h of rs] (rn) {$\ell_r$};
\node[rounded corners,draw,fit=(rs.north)(rn)] {};

\draw[color=blue, dashed,arc] (a1.center) -- (c2.center);
\draw[color=blue, dashed,arc] (c2.center) -- (e3.center);
\draw[color=blue, dashed,arc] (e3.center) -- (h3.center);
\draw[color=blue, dashed,arc] (h3.center) -- (i0.center);
\draw[color=blue, dashed,arc] (i0.center) -- (k1.center);
\draw[color=blue, dashed,arc] (k1.center) -- (n2.center);
\draw[color=blue, dashed,arc] (n2.center) -- (o0.center);
\draw[color=blue, dashed,arc] (o0.center) -- (q1.center);
\draw[color=blue, dashed,] (q1.center) -- ++ (0.5*\l,0);
\draw[color=blue, dashed,<-,>=latex] (a1.center) -- ++ (-0.5*\l,0);

\draw[color= red, arc, dotted, thick] (b0.center) -- (c0.center);
\draw[color= red, arc, dotted, thick] (c0.center) -- (e1.center);
\draw[color= red, arc, dotted, thick] (e1.center) -- (g1.center);
\draw[color= red, arc, dotted, thick] (g1.center) -- (i2.center);
\draw[color= red, arc, dotted, thick] (i2.center) -- (k3.center);

\draw[arc] (f0.center) -- (h0.center);
\draw[arc] (f1.center) -- (h1.center);
\draw[arc] (f2.center) -- (h2.center);
\draw[arc] (f3.center) -- (h3.center);

\draw[arc] (j0.center) -- (l0.center);
\draw[arc] (j1.center) -- (l0.center);
\draw[arc] (j2.center) -- (l0.center);
\draw[arc] (j3.center) -- (l0.center);

\draw[arc] (m0.center) -- (o1.center);
\draw[arc] (m1.center) -- (o2.center);
\draw[arc] (m2.center) -- (o3.center);

\end{tikzpicture}
\end{center}
\caption{Illustration of the digraph $D$ of Section~\ref{subsec:ip-ar} for $\maintenance=4$. Only a few arcs of $D$ (plain arcs) are displayed. The dashed cycle provides a route satisfying the maintenance requirement. The dotted path cannot be completed to a cycle since there is no outgoing arc from $(\ell_k,4)$: It implies that there is no route of the form $\ell_b,\ell_c,\ell_e,\ell_g,\ell_i,\ell_k,\ldots$ that satisfies the maintenance requirement.}
\label{fig:ARgraph}
\end{figure}


We first explain how the problem can be modeled as a disjoint cycle problem in a directed graph. This will make the description of the integer program a straightforward task.

Define the directed graph $D=(V,A)$ as follows. Its vertex set is $\L\times[\maintenance]$. In other words, each flight leg is duplicated $\maintenance$ times.
Each vertex $(\ell,\delta)$ corresponds to a flight leg $\ell\in\L$ with the number of days $\delta\in[\maintenance]$ since the last night spent in a base.
An ordered pair $\left((\ell,\delta),(\ell',\delta')\right)$ is in $A$ if
$(\ell, \ell')$ is an airplane connection  
and we are in one of the three following situations:
\begin{itemize}
\item $\ell$ and $\ell'$ are performed during a same day and $\delta=\delta'$, as illustrated between legs $\ell_f$ and $\ell_h$ on Figure~\ref{fig:ARgraph},
\item $\ell$ and $\ell'$ are not performed on the same day, the airport is a base, and $\delta'=1$, as illustrated between legs $\ell_j$ and $\ell_l$,
\item $\ell$ and $\ell'$ are not performed on the same day, the airport is not a base, and $\delta' - \delta\geq 0$ is the number of days between the arrival of $\ell$ and the departure of $\ell'$, as illustrated between $\ell_m$ and $\ell_o$.  
\end{itemize}
In other words, an arc corresponds to two flight legs that can be consecutive in a route, with the suitable restrictions on the number of days since the last night spent in a base. 
A cyclic sequence of legs $\ell_1,\ldots,\ell_k$ satisfies the maintenance requirement if and only if there exist $\delta_i$ for $i$ in $\{1,\ldots,k\}$ such that $(\ell_1,\delta_1),\ldots,( \ell_k,\delta_k)$ is a cycle in $D$. 
Indeed, suppose that a route $\ell_1,\ldots,\ell_k$ satisfies the maintenance requirement, and denote by $\delta_i$ the number of days since the last night spent in a base before $\ell_i$. Then the definition of $D$ ensures that $(\ell_1,\delta_1),\ldots,( \ell_k,\delta_k)$ is a cycle in $D$. 
Conversely, suppose that a route does not satisfies the maintenance requirement. It contains a sequence of legs $\ell_1,\ldots,\ell_k$ spending more than $\maintenance$ successive days out of a base. Let $\ell_j$ be the last leg in that sequence before the $\maintenance$th night.
Then any path $(\ell_1,\delta_1),\ldots,(\ell_j,\delta_j)$ in $D$ necessarily ends in vertex $(\ell_j,\maintenance)$, which has no outgoing arc, and there is no cycle in $D$ corresponding  to the route. 

Figure~\ref{fig:ARgraph} illustrates such a directed graph $D$. The route
$\ell_a, \ell_c, \ell_e, \ell_h, \ell_i, \ell_k, \ell_n, \ell_o, \ell_q$
satisfies the maintenance requirement and corresponds to the dashed cycle in $D$.
The route
$\ell_b, \ell_c, \ell_e, \ell_g, \ell_i, \ell_k,\ell_n,\ell_o,\ell_q$ does not satisfy the maintenance requirement. The dotted path is an attempt to make it a cycle in $D$, but since there is no arc outgoing from $(\ell_k,\maintenance)$, it is not possible.
 
We choose arbitrarily one instant in the week and we denote by $A_0$ the set of arcs $((\ell,\delta), (\ell',\delta'))$ ``crossing this instant'', i.e., such that the time interval between the departure of $\ell$ (included) and the departure of $\ell'$ (excluded) contains the instant.
 Define moreover $V_{\ell}$ to be the set $\big\{(\ell,\delta)\in V\colon \delta\in[\maintenance]\big\}$.

\begin{prop}\label{prop:AR}
Feasible solutions of the aircraft routing problem are in one-to-one correspondence with collections $\C$ of vertex disjoint cycles in $D$ such that we have simultaneously
\begin{enumerate}[label=\textup{(\roman*)}]
\item \label{intersec} for each $\ell$, exactly one cycle in $\C$ has a nonempty intersection with $V_{\ell}$, and this intersection consists of a single vertex.
\item \label{A0} $\C$ has at most  $n^{\mathrm{a}}$ arcs in $A_0$. 
\end{enumerate}
\end{prop}


\begin{figure}
	\center
	\begin{tikzpicture}

\tikzset{vert/.style={draw, circle, inner sep=0pt, text width=1mm}}
\tikzset{arc/.style={->}}

\draw [blue, fill= blue!10] plot [smooth cycle] coordinates {(-1.9,1.4) (0,1.6) (1.9,1.4) (2.5,0) (1.9,-1.4) (0, -1.6) (-1.9, -1.4) (-2.5,0)};
\draw [blue, fill=white] plot [smooth cycle] coordinates {(-1.3,0.6) (0,0.8) (1.3,0.6) (1.6,0) (1.3,-0.6) (0, -0.8) (-1.3, -0.6) (-1.6,0)};

\draw[thick] (0, 0.7) -- (0,1.7);
\draw[dotted] (1.3, 0.5) -- (2.1, 1.4);
\draw[dotted] (1.45,-0.35) -- (2.5,-0.45);
\draw[dotted] (0.8,-0.65) -- (0.9,-1.65);
\draw[dotted] (-0.8,-0.65) -- (-0.9,-1.65);
\draw[dotted] (-1.45,-0.35) -- (-2.5,-0.45);
\draw[dotted] (-1.3, 0.5) -- (-2.1, 1.4);

\node[vert] (v1) at (0.2,1.4) {};
\node[vert] (v2) at (1.7,1.4) {};
\node[vert] (v3) at (1.7,-0.9) {};
\node[vert] (v4) at (-1.6,-1.2) {};
\node[vert] (v5) at (-1.7,0.6) {};
\node[vert] (v6) at (-1,0.9) {};
\node[vert] (v7) at (-0.5,0.9) {};
\node[vert] (v8) at (1.3,0.9) {};
\node[vert] (v9) at (2.2,0) {};
\node[vert] (v10) at (2,-1.2) {};
\node[vert] (v11) at (0,-1.4) {};
\node[vert] (v12) at (-2.4,0) {};
\node[vert] (v13) at (-2,1.1) {};

\node[above right=-0.1cm and -0.2cm of v7, red] {$A_0$};

\draw[arc] (v1) -- (v2);
\draw[arc] (v2) -- (v3);
\draw[arc] (v3) -- (v4);
\draw[arc] (v4) -- (v5);
\draw[arc] (v5) -- (v6);
\draw[arc] (v6) -- (v7);
\draw[arc, dashed, red] (v7) -- (v8);
\draw[arc] (v8) -- (v9);
\draw[arc] (v9) -- (v10);
\draw[arc] (v10) -- (v11);
\draw[arc] (v11) -- (v12);
\draw[arc] (v12) -- (v13);
\draw[arc, dashed, red] (v13) -- (v1);



\end{tikzpicture}
	\caption{A two-week route crosses $A_0$ twice}
	\label{fig:timeline}
\end{figure}



Therefore, the aircraft routing problem is equivalent to deciding whether the following integer program has a feasible solution:

\begin{subequations}
        \makeatletter
        \def\@currentlabel{AR}
        \makeatother
\label{eq:ARprogram}
        \renewcommand{\theequation}{AR.\arabic{equation}}
\begin{alignat}{2}
&\ds{\sum_{a\in\delta^-(v)}x_a = \sum_{a\in\delta^+(v)}x_a} && \qquad\forall v\in V \smallskip \label{eq:arflow}\\
&\ds{\sum_{a\in\delta^-(V_{\ell})}x_a = 1} && \qquad\forall \ell\in\L \smallskip \label{eq:partition} \\
&\ds{\sum_{a\in A_{0}}x_a \leq n^{\mathrm{a}}} && \smallskip \label{eq:number_airplanes}\\
&x_a \in \{0,1\} &&\qquad \forall a\in A. \label{eq:ar01}
\end{alignat}
\end{subequations}
Equation~\eqref{eq:arflow} is the flow equation. Together with \eqref{eq:partition}, it ensures that the solution is composed of vertex disjoint cycles. Equation~\eqref{eq:partition} ensures that \ref{intersec} is satisfied and Equation~\eqref{eq:number_airplanes} ensures that \ref{A0} is satisfied.

\subsection{Solution method}

{
The solution we propose is to implement directly the integer program \eqref{eq:ARprogram} in any standard MIP solver. Its number of constraints is $|\L|\times\maintenance+|\L|+1$ and its number of variables is the number of airplane connections times $\maintenance$. Program \eqref{eq:ARprogram} is therefore of tractable size, and current off-the-shelf solvers can solve industrial instances in a few minutes.
See our experiments in Section~\ref{sec:experimentalResults}.}

\section{Column generation approach to crew pairing}
\label{sec:crewPairing}

\subsection{Problem formulation} 
\label{sub:problem_formulation}

Roughly speaking, the crew pairing problem is similar to the aircraft routing problem: instead of building sequences of flight legs for the airplanes (the routes), crew pairing requires to build {\em pairings}, which are sequence of flight legs operated by the crews.
While routes are cyclic sequences, pairings are noncyclic sequences (and they are often quite short). The set of flight legs have to be partitioned into pairings, but the constraints are much more complicated. Before stating formally the crew pairing problem, we introduce some terminology.

A pair of flight legs $(\ell,\ell')$ is a {\em connection} if it satisfies:
\begin{itemize}
\item The arrival airport of $\ell$ is the departure airport of $\ell'$
\item The departure time of $\ell'$ minus the arrival time of $\ell$ is bounded from below by a fixed quantity (which can depend on the airport, the time, and the fleet). This quantity is in general different from the similar one for aircraft routing.
\end{itemize} If the arrival of $\ell$ and the departure of $\leg'$ are on the same day, then it is a \emph{day connection}. Otherwise, it is a \emph{night connection}. If the duration of a night connection is smaller than a threshold, then it is a \emph{reduced rest}. (This term is due to the rest taken by crews performing a night connection.)

A {\em pairing} is a sequence of distinct flight legs such that any two consecutive flight legs form a connection. The subsequence of a pairing formed by all flight legs operated during a same day is a {\em duty}.

To be {\em feasible}, a pairing $\ell_1,\ldots,\ell_k$ has to satisfy the following rules:
\begin{enumerate}[label=(\alph*)]
	\item the period between the departure of $\ell_1$ and the arrival of $\ell_k$ spans at most $4$ days, \label{rule:length}
	\item $\leg_{1}$ starts and $\ell_{k}$ ends in one of the Paris airports, \label{rule:Paris}
	\item each duty contains at most $4$ flight legs. 
If a duty starts with a leg $\ell'$, and the night connection $(\ell,\ell')$ that leads to $\ell'$ is a reduced rest,
	then the number of legs of the duty is at most $3$, \label{rule:NbLegs}
	\item the total flying duration in a duty does not exceed $F(t)$, where $F$ is a given function and $t$ is the time at which the first leg of the duty departs, \label{rule:FlyingTime}
\end{enumerate}
as well as more than $70$ other rules which encode the IR-OPS regulation of the European Aviation Safety Agency  and Air France working rules. A pairing is \emph{long} if it spans $4$ days. A duty is \emph{long} if it contains more than $3$ flight legs, and \emph{short} otherwise. We denote by $\pairingSet$ the set of feasible pairings

Operating a pairing $p$ in $\pairingSet$ has a cost $c_{p}$ that corresponds to crew wages and hotel nights. Given a set $\legSet$ of flight legs, solving the crew pairing problem consists in selecting a collection of feasible pairings of minimum total cost so that each leg $\leg$ in $\legSet$ belongs to exactly one of them, and so that the following global constraints are satisfied: The proportion of long pairings in the solution is less than or equal to a quantity $\alpha$, and the proportion of long duties is at most a quantity $\beta$.

Our purpose here is to introduce the main modeling ideas and not to get into details of the intricacies of the IR-OPS and Air France regulations. In the rest of the paper, we therefore present these ideas on a simplified problem with only the four illustrating rules \ref{rule:length}, \ref{rule:Paris}, \ref{rule:NbLegs}, and \ref{rule:FlyingTime}.
All the other rules are also taken into account in the numerical results.

\subsection{Modeling as an integer program}


As pairings must satisfy many non-linear rules such as rule \ref{rule:NbLegs},
it is difficult to model crew pairing using a compact integer program that has both a good linear relaxation and a tractable size.
The literature therefore generally uses a column generation approach where rules complexity are hidden in the set of variables (and dealt with using ad-hoc algorithms in the pricing subproblem). We also use a column generation approach.

 The binary variable $y_{p}$ indicates if a pairing $p$ in $\pairingSet$ belongs to the solution. 

\begin{equation}\tag{CP}\label{eq:MasterProblem}
	\begin{array}{rl@{\qquad}l}
	\min & \displaystyle \sum_{\pairing \in \pairingSet} c_{\pairing}y_{\pairing} \\
	\mathrm{s.t.}& \displaystyle \sum_{ \pairing \ni \leg} y_{\pairing} = 1 & \forall \leg \in \legSet\\
	& \displaystyle \sum_{\pairing \in \pairingSet^{\mathrm{l}}}y_{\pairing} \leq  \lpProp\sum_{\pairing\in\pairingSet} y_{\pairing}   \\
	& \displaystyle \sum_{\pairing \in \pairingSet} \left((1 - \ldProp ) \Delta^{\mathrm{l}}(p) - \ldProp \Delta^{\mathrm{s}}(p)\right) y_{\pairing} \leq 0\\
	& y_{\pairing}\in \{0,1\} & \forall \pairing \in \pairingSet,
	\end{array}
\end{equation}
where $p\ni\ell$ means that the flight leg $\ell$ is present in $p$, where $\pairingSet^{\mathrm{l}}$ is the set of long pairings, and where $\Delta^{\mathrm{s}}(p)$ (resp.~$\Delta^{\mathrm{l}}(p)$) is the number of short (resp.~long) duties in a pairing $p$. The first constraint ensures that each leg is covered, the second that the proportion of long pairings is less than or equal to $\alpha$, and the third that the proportion of long duties is less than or equal to~$\beta$. 

\subsection{Column generation approach} 
\label{sub:column_generation_approach}

We propose an exact method for solving the program~\eqref{eq:MasterProblem}. It is based on column generation. 
We describe the method without assuming special knowledge in column generation. We will in particular be sketchy on the theoretical rationale; more details on that topic can be found for instance in a survey by \citet{lubbecke2011column}.





\begin{algorithm}
\begin{algorithmic}[1]
\STATE \textbf{initialize} $\mathcal{P}'$ in such a way that \eqref{eq:MasterProblem} restricted to $\mathcal{P}'$ is feasible (e.g., taking all possible pairings of two flight legs makes the job);
\REPEAT 
\STATE solve the linear relaxation of \eqref{eq:MasterProblem} restricted to $\mathcal{P}'$ with any standard solver; 	\label{step:master}
\STATE denote by $c^{\mathrm{low}}$ its optimal value;
\STATE find a pairing $p$ of minimum reduced cost $\tilde{c}_{p}$; {\em (pricing subproblem)}\label{step:minreduced}
\IF{($\tilde{c}_{p}<0$)}
\STATE add $p$ to $\mathcal{P}'$;
\ENDIF
\UNTIL ($\tilde{c}_{p}\geq 0$ for all $p\in\mathcal{P}$)
\STATE solve \eqref{eq:MasterProblem} restricted to $\pairingSet'$ with any standard solver; \label{step:firstUB}
\STATE denote by $c^{\mathrm{upp}}$ its optimal value; \label{step:firstInteger}
\STATE add to $\pairingSet'$ all pairings with reduced cost (from the last linear program of Step~\ref{step:master}) non-larger than $c^{\mathrm{upp}} - c^{\mathrm{low}}$; \label{step:completeWithPathGap} 
\STATE solve \eqref{eq:MasterProblem} restricted to $\pairingSet'$ with any standard solver;\label{step:final}
\STATE \textbf{return} its optimal solution~$\yy^*$; 
\end{algorithmic}
\caption{Column generation algorithm}
\label{alg:CG}
\end{algorithm}



Algorithm~\ref{alg:CG} describes our column generation approach, which maintains a subset of pairings $\mathcal{P}' \subseteq \mathcal{P}$. The idea of column generation is to solve the {\em master problem}, which is the linear relaxation of \eqref{eq:MasterProblem} on such a subset $\mathcal{P}'$, and to check if by chance the optimal solution found on this restricted version is also an optimal solution of the master problem with the full set $\mathcal{P}$. This checking is done exactly as in the classical simplex algorithm: It is the optimal solution for the full problem if all reduced costs are nonnegative. Since the number of elements in $\mathcal{P}$ is huge, it is not possible to compute and check all these reduced costs one by one. However, given an element $p$ in $\mathcal{P}$, it is always possible to compute its reduced cost from the value of the dual variables by standard linear programming theory. To find the element $p$ with the smallest reduced cost, an auxiliary optimization problem instantiated by the values of the dual variables is solved: The {\em pricing subproblem}. In Algorithm~\ref{alg:CG}, this is done in Step~\ref{step:minreduced}. The exact method to solve the pricing subproblem is described in Section~\ref{sec:pricing_subproblem}. Right before Step~\ref{step:firstUB}, the linear relaxation of \eqref{eq:MasterProblem} is fully solved.

Since the total number of possible pairings is finite, Step~\ref{step:master} -- which consists in solving the master problem -- is repeated only finitely many times, and thus the overall method terminates in finite time. After having performed Step~\ref{step:master} for the last time, $c^{\mathrm{low}}$ is a lower bound on the optimal value of \eqref{eq:MasterProblem}. Step~\ref{step:firstUB} provides a first feasible solution of \eqref{eq:MasterProblem}. The value $c^{\mathrm{upp}}$ (Step~\ref{step:firstInteger}) is thus an upper bound on the optimal value of \eqref{eq:MasterProblem}. At that time, we have thus a lower bound on the optimal value, and a feasible solution.

The purpose of the remaining steps is to ``close the gap''. The idea consists in generating all pairings $p$ that might be in an optimal solution. These pairings are precisely those whose reduced cost is smaller than or equal to $c^{\mathrm{upp}} - c^{\mathrm{low}}$. This is completely formalized by Lemma~\ref{lem:baldacci} below. Finding all these pairings in Step~\ref{step:completeWithPathGap} is a variant of the pricing subproblem, briefly discussed in Remark~\ref{rem:modifiedSubproblem} of Section~\ref{sec:pricing_subproblem}.
 At the end, the solution $\yy^*$ is an optimal solution of \eqref{eq:MasterProblem}. 

\begin{lem}[{\citet[Proposition 2.1, p.~389]{nemhauser1988integer}}]
\label{lem:baldacci}
Consider an integer program in standard form with variables $(z_i)$ for which the linear relaxation admits a finite optimal value $\bar v$. Suppose given an upper bound $\UB$ on the optimal value of the integer program. Then for every $i$ such that $\tilde{c}_i>\UB-\bar v$, the variable $z_i$ is equal to $0$ in all optimal solutions of the integer program, where $\tilde{c}_i$ denotes the reduced cost of the variable $z_i$ when the linear relaxation has been solved to optimality. 
\end{lem}

\begin{rem}
The last steps of Algorithm~\ref{alg:CG} works only if the gap $c^{\mathrm{upp}} - c^{\mathrm{low}}$ is small, as otherwise a huge number of pairings might be added at Step~\ref{step:completeWithPathGap}. We found out numerically that it is the case for all Air France instances, and we thus use this technique. 
Lemma~\ref{lem:baldacci} has been recently used by \citet{cacchiani2016optimal} on the integrated problem. They underline that when $c^{\mathrm{upp}} - c^{\mathrm{low}}$ is large, a branch and bound approach is required.
\end{rem}

\section{Pricing subproblem} 
\label{sec:pricing_subproblem}

We now introduce a solution scheme for the pricing subproblem
\begin{equation}\label{pricing}
	\min_{\pairing \in \pairingSet} \tilde{c}_{p}. 
\end{equation}
As pairings can be considered as paths satisfying constraints in 
the graph whose vertices are the legs and arcs the connections, the pricing subproblem is generally solved as a resource constrained shortest path problem, and we do not depart from this approach. We model it within the \MRCSP framework \citep{parmentier2015algorithms}, which we now briefly describe. This framework is rather abstract, but practically, it only requires to implement a few operators on the resource set.
This work is the first application of the \MRCSP framework to an industrial problem.

\subsection{Framework and algorithm} 
\label{sub:framework_and_algorithm}


A binary operation $\rplus$ on a set $\rset$ is \emph{associative} if $\re \rplus (\re' \rplus \re'') = (\re \rplus \re')\rplus \re''$ for $\re,\re',$ and $\re''$ in $\rset$. An element $0$ is \emph{neutral} if $0 \rplus \re =  \re \rplus 0 = \re$ for any $\re$ in $\rset$. A set $(\rset,\rplus)$ is a \emph{monoid} if $\rplus$ is associative and admits a neutral element. A partial order $\rleq$ is \emph{compatible} with $\rplus$ if the mappings $\re \mapsto \re \rplus \re'$ and $\re \mapsto \re' \rplus \re$ are non-decreasing according to this order for all $\re' $ in $\rset$. A partially ordered set $(\rset,\rleq)$ is a \emph{lattice} if any pair $(q,q')$ of elements of $\rset$ admits a greatest lower bound or \emph{meet} denoted by $\re\meet\re'$, and a least upper bound or \emph{join} denoted by $\re\join\re'$. A set $(\rset,\rplus,\rleq)$ is a \emph{lattice ordered monoid} if $(\rset,\rplus)$ is a monoid, $(\rset,\rleq)$ is a lattice, and $\rleq$ is compatible with $\rplus$.

Given a digraph $D=(V,A)$, a lattice ordered monoid $(\rset,\rplus,\rleq)$, elements $q_a\in\rset$ for each $a\in A$, origin and destination vertices $o$ and $d$, and two non-decreasing mappings $\rcost : \rset \rightarrow \R$ and $\rmeas:\rset \rightarrow \{0,1\}$, the \MRCSP seeks
$$\text{an $o$-$d$ path $P$ of minimum $\rcost\left(\bigrplus_{a\in P}\re_{a}\right)$ among those satisfying $\rmeas\left(\bigrplus_{a\in P}\re_{a}\right) = 0$,} $$
where $\bigrplus_{a\in P}$ is always performed in the order of the arcs on the path $P$ (the operation $\oplus$ is not necessarily commutative).
We call such a $q_a$ the {\em resource} of the arc $a$.
The sum $\bigrplus_{a\in P}\re_{a}$ is the \emph{resource} of a path $P$, and we denote it by $\re_{P}$. The real number $\rcost\left(\re_P\right)$ is its \emph{cost}, and the path $P$ is \emph{feasible} if $\rmeas\left(\re_P\right)$ is equal to $0$. We therefore call $\rcost$ and $\rmeas$ the {\em cost} and the {\em infeasibility functions}. 
\def\key{\operatorname{key}}

We now describe an \emph{enumeration algorithm} for the \MRCSP. It follows the standard labeling scheme \cite{irnich2005shortest} for resource constrained shortest paths. 
The specificity of our algorithm is that it uses, for each $v$ in $V$, a set $B_{v}$ of bounds such that, 
\begin{equation}\label{eq:boundSetProp}
	\text{for each $v$-$d$ path $Q$, there is a $b\in B_{v}$ with $b\rleq \re_{Q}$.}
\end{equation}
The lattice ordered monoid framework enables to design procedures to build these sets of bounds; see \Cref{sub:enhanced_bounds_computations_for_the_column_generation_context} for more details. 
Having defined these bounds, we define $\key(P)$ as 
\begin{equation}\label{eq:keyDefinition}
	\key(P) = \min \{c(\re_P \rplus b)\colon b \in B_v, \rmeas(\re_P \rplus b)=0\} \quad \text{where $v$ is the last vertex of $P$.}
\end{equation}
The empty path at a vertex $v$ is the path with no arcs starting and ending at vertex $v$.
By definition of paths resources, its resource is the neutral element of the monoid.
A path $P$ {\em dominates} a path $Q$ if $\re_P\rleq\re_Q$.
During the algorithm, a list $\mathsf{L}$ of partial paths, an upper bound $c_{od}^{UB}$ on the cost of an optimal solution, and
lists $(\mathsf{L}_{v}^{\mathrm{nd}})_{v \in V}$ of non-dominated $o$-$v$ paths 
 are maintained. Algorithm~\ref{alg:enumeration} states our algorithm. We denote by $P+a$ the path composed of a path $P$ followed by an arc $a$.

\begin{algorithm}
\begin{algorithmic}[1]
\STATE \textbf{input}: sets $(B_v)_{v \in V}$ satisfying~\eqref{eq:boundSetProp}; \quad \emph{(see Section~\ref{sub:enhanced_bounds_computations_for_the_column_generation_context})}
\STATE \textbf{initialization:}  $c_{od}^{UB} \leftarrow +\infty$, $\mathsf{L} \leftarrow \emptyset$, and $\mathsf{L}_{v}^{\mathrm{nd}} \leftarrow \emptyset$ for each $v \in V$;
\STATE add the empty path at the origin $o$ to $\mathsf{L}$ and $\mathsf{L}_{o}^{\mathrm{nd}}$;
\WHILE{$\mathsf{L}$ is not empty}
\STATE $P \leftarrow$ a path of minimum $\key(P)$ in $\mathsf{L}$; \label{step:pathSelection} 
\STATE $\mathsf{L} \leftarrow \mathsf{L}\backslash \{P\}$;
\STATE $v\leftarrow$ last vertex of $P$;
\IF{$v=d$, $\rmeas(\re_P) = 0$, and $\rcost(\re_P)<c_{od}^{UB}$} \label{step:extStart}
\STATE $c_{od}^{UB} \leftarrow \rcost(\re_{P}) $; \label{step:codub}
\ELSE \quad\quad\quad\quad\emph{(extension of $P$)}
\FORALL{$a\in\delta^+(v)$}\label{step:arcLoop}
\STATE $Q \leftarrow P+a$; \label{step:Qdef}
\STATE $w\leftarrow $ last vertex of $Q$; \label{step:Qstarts}
\IF{$\exists b \in B_w$ such that $\rmeas(\re_{Q} \rplus b) = 0$ and $\rcost(\re_{Q} \rplus b) < c_{od}^{UB}$} \label{step:LBtest}
\IF{$Q$ is not dominated by any path in $\mathsf{L}_{w}^{\mathrm{nd}}$} \label{step:DomTest} 
\STATE $\mathsf{L}_{w}^{\mathrm{nd}} \leftarrow \mathsf{L}_{w}^{\mathrm{nd}} \cup \{Q\}$ and remove from $\mathsf{L}_{w}^{\mathrm{nd}}$ and $\mathsf{L}$ every path dominated by $Q$; \label{step:removal}
\STATE $\mathsf{L} \leftarrow \mathsf{L} \cup \{Q\}$; \label{step:updateL}
\ENDIF 
\ENDIF \label{step:Qends}
\ENDFOR
\ENDIF \label{step:extEnd}
\ENDWHILE
\STATE \textbf{return} $c_{od}^{UB}$;
\end{algorithmic}
\caption{Enumeration algorithm for the \MRCSP}
\label{alg:enumeration}
\end{algorithm}

\begin{prop}\label{prop:MRCSPalgoConvergence}
Suppose that $D$ is acyclic. Then Algorithm~\ref{alg:enumeration} converges after a finite number of iterations, and, at the end of the algorithm, $c_{od}^{UB}$ is equal to the cost of an optimal solution of the \MRCSP if such a solution exists, and to $+\infty$ otherwise. 
\end{prop}

The specificity of our approach lies in the use of the bounds of Equation~\eqref{eq:boundSetProp} in the algorithm. While it is well known that the use of lower bounds is a key element in the performance of the enumeration algorithms \cite{dumitrescu2003improved}, our approach is the first to allow the use of bounds with non-linear constraints such as~\ref{rule:NbLegs} and~\ref{rule:FlyingTime}. Not only these bounds are used to discard more paths, but they are also used to improve the order in which the paths are considered by the algorithm. The two main resource constrained shortest path algorithms in the literature \citep{irnich2005shortest} differ by the order in which they consider paths. 
The \emph{label correcting} algorithm is obtained from our one by using $c(\re_P)$ as $\key(P)$ and removing the test of Step~\ref{step:LBtest}. It is for instance described by \citet{dunbar2012robust,dunbar2014integrated} in the context of crew pairing. The \emph{label setting} algorithm considers vertices $v$ in a topological order, and then apply Steps~\ref{step:extStart} to \ref{step:extEnd} for each path $P$ in $\mathsf{L}_v^{\mathrm{nd}}$. Again, the test of Step~\ref{step:LBtest} is not present. 

Finally, as is has already been noted, we go further by using sets of bounds rather than singletons (we give additional explanations in Section~\ref{sub:enhanced_bounds_computations_for_the_column_generation_context}).%
\begin{rem}\label{rem:modifiedSubproblem}
Step~\ref{step:completeWithPathGap} of Algorithm~\ref{alg:CG} requires to solve the following variant of the \MRCSP:
 $$\text{generate all the $o$-$d$ paths $P$ satisfying $\rmeas(\re_{P}) = 0$ and $\rcost(\re_{P}) \leq c^{\mathrm{upp}} - c^{\mathrm{low}}$.}$$ 

Algorithm~\ref{alg:enumeration} can be easily adapted to this variant. It suffices to maintain a set $\mathsf{S}$ of solutions (initially empty), to replace $c_{od}^{UB}$ by $c^{\mathrm{upp}} - c^{\mathrm{low}}$ in Steps~\ref{step:extStart} and~\ref{step:LBtest}, to replace Step~\ref{step:codub} by $\mathsf{S} \leftarrow \mathsf{S} \cup \{P\}$, and to return $\mathsf{S}$. The set $\mathsf{S}$ returned contains all the $o$-$d$ paths $P$ satisfying $\rmeas(\re_{P}) = 0$ and $\rcost(\re_{P}) \leq c^{\mathrm{upp}} - c^{\mathrm{low}}$.
\end{rem}

\begin{rem}
The terminology ``label setting'' and ``label correcting'' varies in the literature. We stick here to \citet{irnich2005shortest}. What \citet{dunbar2012robust,dunbar2014integrated} call a ``label setting'' algorithm is a label correcting algorithm according to \citet{irnich2005shortest}.
\end{rem}


\subsection{Modeling the pricing subproblem} 
\label{sub:modeling_the_pricing_subproblem}

We now explain how to model our pricing subproblem~\eqref{pricing} in the \MRCSP framework.
As already mentioned, we only consider rules \ref{rule:length} to \ref{rule:FlyingTime} (Section~\ref{sub:problem_formulation}) to focus on ideas rather than on the full intricacies of the regulation. In other words, the set of feasible pairings is the set of pairings satisfying  rules \ref{rule:length} to \ref{rule:FlyingTime}.
For simplicity, we also omit in the master problem the long pairings and long duties constraints, and we assume that the cost $c_p$ is of the form $\sum_{(\ell,\ell') \in p}c_{(\ell,\ell')}$. 
We emphasize that all IROPS and Air France rules, as well as the real costs and the long pairings and long duties constraints, are taken into account in the numerical experiments.
The reduced cost is then of the form $\tilde{c}_{p} = \sum_{(\ell,\ell') \in p}c_{(\ell,\ell')} +  \sum_{\ell \in p}z_{\ell}$, where $z_\ell$ is the dual variable associated to the partitioning constraint. 

We now model this toy subproblem as a \MRCSP. According to Section~\ref{sub:framework_and_algorithm}, we have to describe the digraph with its origin and destination vertices, the lattice ordered monoid, the resources on the arcs, and the cost and infeasibility functions. We actually solve a shortest path problem for each sequence of four consecutive days in a week, in order to satisfy rules~\ref{rule:length} and~\ref{rule:Paris} (Section~\ref{sub:problem_formulation}). We thus solve seven \MRCSP instances per iteration of the pricing subproblem.


\subsubsection{The digraph with its origin and destination vertices} Let $D = (V,A)$ be the acyclic digraph defined as follows. The vertex set is $V= \legSet \cup \{o,d\}$, where $\legSet$ is the set of legs of four consecutive days, $o$ is a dummy origin vertex, and $d$ is a dummy destination vertex. The arc set $A$ contains an arc $(o,\ell)$ for all legs $\ell$ starting in Paris on the first day of the period, an arc $(\ell,d)$ for all legs $\ell$ ending in Paris, and an arc $(\ell,\ell')$ for each connection $(\ell,\ell')$. With these definitions, pairings starting on the first day of the period and satisfying rules~\ref{rule:length} and \ref{rule:Paris} are in one-to-one correspondence with $o$-$d$ paths in $D$.

\subsubsection{The lattice ordered monoid}
The monoid $\rset$ we use for the resources is of the form  
$\rset^{\rmeas} \times \R$, where $ \rset^{\rmeas} = (\Z_{+} \times \R_{+}) \cup (\Z_{+} \times \R_{+})^{2} \cup \{\infty\}$. 

An element $(n,f)\in\Z_{+} \times \R_{+}$ models the resource used over a single day by a pairing: $n$ flight legs and a total flying duration $f$. An element $(n^b,f^b,n^e,f^e)\in(\Z_{+} \times \R_{+})^{2}$ models the resource used over the first and last days of a pairing lasting more than one day.
Quantities~$n^b$ and $f^b$ are the number of flight legs and flying duration on the first day of a pairing, and $n^e$ and $f^e$ are the number of flight legs and flying duration on the last day of a pairing.
  The element $\infty$ is used to capture infeasibility of certain pairings.

Let $F_{\mathrm{m}} = \displaystyle\max_t F(t)$. We define the operator $\rplus$ on $\rset^{\rmeas}$ as follows.
 \begin{align*}
	r \rplus \infty &= \infty \rplus r = \infty \quad \text{for all } r \in\rset^\rmeas\\
	(n,f) \rplus (\tilde{n},\tilde{f}) &= (n + \tilde{n}, f + \tilde{f})\\
	(n,f) \rplus (\tilde{n}^{b},\tilde{f}^{b},\tilde{n}^{e},\tilde{f}^{e}) &= (n + \tilde{n}^{b}, f + \tilde{f}^{b}, \tilde{n}^{e}, \tilde{f}^{e}) \\
	(n^{b},f^{b},n^{e},f^{e}) \rplus (\tilde{n},\tilde{f}) &= (n^{b}, f^{b}, n^{e} + \tilde{n}, f^{e}+\tilde{f})\\
	(n^{b},f^{b},n^{e},f^{e}) \rplus (\tilde{n}^{b},\tilde{f}^{b},\tilde{n}^{e},\tilde{f}^{e}) &=  \left\{
	\begin{array}{ll}
	 \infty & \text{if $n^{e} + \tilde{n}^{b} >4$ or $f^{e} + \tilde{f}^{b}>F_{\mathrm{m}}$,} \\
	 (n^{b},f^{b}, \tilde{n}^{e}, \tilde{f}^{e}) & \text{otherwise.}
	 \end{array} \right.
\end{align*}




We define $\rleq$ on $\rset^{\rmeas}$ by
\begin{align*}
(0,0)\rleq \re \quad &\text{and} \quad \re \rleq \infty \quad \text{for all } r \in \rset^\rmeas \\
(n,f) \rleq (\tilde{n},\tilde{f}) \quad &\text{if} \quad n\leq \tilde{n} \quad \text{and} \quad f \leq \tilde{f} \\
(n^{b},f^{b},n^{e},f^{e}) \rleq (\tilde{n}^{b},\tilde{f}^{b},\tilde{n}^{e},\tilde{f}^{e}) \quad & \text {if} \quad n^{b} \leq \tilde{n}^{b}, \enskip f^{b} \leq \tilde{f}^{b},  \enskip n^{e} \leq \tilde{n}^{e}, \text{ and } f^{e} \leq \tilde{f}^{e},
\end{align*}
and a pair $(n,f) \neq (0,0)$ is not comparable with $(n^{b},f^{b},n^{e},f^{e})$.


\begin{lem}\label{lem:monoid}
$(\rset^{\rmeas}, \rplus, \rleq)$ is a lattice ordered monoid.
\end{lem}

As $(\R,+,\leq)$ is a lattice ordered monoid, the monoid \emph{$\rset=\rset^{\rmeas}\times \R$ is a lattice ordered monoid when endowed with the componentwise sum and order.}

\subsubsection{Resources on the arcs}
Consider an arc $(\ell,\ell')$ of $D$. If it is a day connection, then we define its resource to be $\left((1,f(\ell')),c_{(\ell,\ell')} + z_{\ell'}\right)$, where $f(\ell')$ is the flying duration of leg~$\ell'$, and $z_{\ell'}$ is the dual variable of the cover constraint associated to $\ell'$ in \eqref{eq:MasterProblem}. If it is a night connection, then we define its resource to be $\left((0,0,n^{e},f^{e}),c_{(\ell,\ell')} + z_{\ell'}\right)$, where $n^{e} =2$ if $(\ell,\ell')$ is a reduced rest, and $1$ otherwise, and $f^{e} = f(\ell') + F_{\mathrm{m}} - F(t)$, where $t$ is the departure time of $\ell'$.
Similarly, each arc $(o,\ell')$ has resource $\big((0,0,1,f^{e}),z_{\ell'}\big)$, and each arc $(\ell,d)$ resource~$\big((0,0,0,0),0\big)$.

%


\subsubsection{Cost and infeasibility functions}
Given $\re = (r,z) \in \rset$, we define
\begin{align*}
\rmeas\left((r,z)\right) = \rmeas_{\rset^{\rmeas}}(r) \quad \text{and} \quad \rcost((r,z)) = z
\end{align*}
\noindent
where $\rmeas_{\rset^{\rmeas}}$ is defined on $\rset^{\rmeas}$ by
\begin{align*}
	\rmeas_{\rset^{\rmeas}}\big((n,f)\big) &= \max\left(\ind_{(4,\infty)}(n), \ind_{(F_{\mathrm{m}},\infty)}(f)\right), \\
	\rmeas_{\rset^{\rmeas}}\left((n^{b},f^{b},n^{e},f^{e})\right) &= \max\left(\ind_{(4,\infty)}(n^{b}), \ind_{(F_{\mathrm{m}},\infty)}(f^{b}),\ind_{(4,\infty)}(n^{e}), \ind_{(F_{\mathrm{m}},\infty)}(f^{e})\right), \\
	\rmeas_{\rset^{\rmeas}}(\infty) &= 1,
\end{align*}
where $\ind_{I}$ denotes the indicator function of a set $I$. With this definition, the feasibility function $\rmeas$ encodes the satisfaction of rules~\ref{rule:NbLegs} and \ref{rule:FlyingTime}.


\subsubsection{Conclusion}
The following proposition concludes the reduction of the pricing subproblem to a \MRCSP.
\begin{prop}\label{prop:pricing}
The sequence of flight legs $p$ corresponding to an $o$-$d$ path $P$ is in $\pairingSet$ (i.e., is a feasible pairing) if and only if $\rmeas(\re_{P}) = 0$. In that case, $\rcost(\re_{P}) = \tilde{c}_{p}$. 
\end{prop}

Appendix~\ref{sec:appendix_example} details how this reduction works on a small pricing subproblem instance and shows a few typical iterations of Algorithm~\ref{alg:enumeration}.

\subsection{Bounds on resources} 
\label{sub:enhanced_bounds_computations_for_the_column_generation_context}

\begin{figure}
a.
\begin{tikzpicture}
\def\l{1.5}
\def\h{1.9}
\node[vert] (o) at (0*\l, 0*\h) {$o$};
\node[vert] (u) at (2*\l, 0*\h) {$w$};
\node[vert] (v) at (3*\l, 1*\h) {};
\node[vert] (w) at (3*\l,-1*\h) {};
\node[vert] (d) at (4*\l, 0*\h) {$d$};
\draw[dashed, ->, gray] (o) to node[midway, above, black] {$(1,1)$} node[midway, below, black] {$Q$}  (u);

\draw[arc] (u) to node[midway, above left] {$(2,0)$}  (v);
\draw[arc] (u) to node[midway, below left] {$(0,2)$}  (w);
\draw[arc] (u) to node[midway, above] {$(2,2)$} node[midway, below] {$R'$}  (d);
\draw[arc] (v) to node[midway, above right] {$(1,1)$}  (d);
\draw[arc] (w) to node[midway, below right] {$(1,1)$}  (d);

\node[below = 0.1 of v] (vs) {$R$};
\draw [dashed, gray,->] plot [smooth, tension=1] coordinates { (u.north east) (vs) (d.north west)};
\node[above = 0.1 of w] (wn) {$R''$};
\draw [dashed, gray,->] plot [smooth, tension=1] coordinates { (u.south east) (wn) (d.south west)};

\end{tikzpicture}
\qquad b.
\begin{tikzpicture}[scale=0.8]
	\def\l{1.5}
	\def\h{1.5}

	\tikzset{axe/.style={->}}
	\draw[axe] (0,0) -- (0,3.6*\h) node[below left]{$\re^{1}$} ;
	\draw[axe] (0,0) -- (3.6*\l,0) node[right]{$\re^{2}$};
	\node (x1) at (3*\l,1*\h) {$\times$};
	\node[right=0 of x1] {$q_{R''}$};
	\node (x2) at (1*\l,3*\h) {$\times$};
	\node[right=0 of x2] {$q_{R}$};
	\node (z) at (2*\l,2*\h) {$\times$};
	\node[right=0 of z] {$q_{R'}$};

	\node[color = blue] (y1) at (1*\l,2*\h) {$\circ$};
	\node[color = blue] (y2) at (3*\l,1*\h) {$\circ$};
	\node[color = red] (m) at (1*\l,1*\h) {$\diamond$};
	\node[below = 0 of m, color=red] {$\re_{R} \meet\re_{R'} \meet \re_{R''}$};
	\node[below left = 0 and -0.6 of y1, color=blue] {$\re_R \meet\re_{R'}$};
	
	\draw[dashed] (z.center) -- (y1.center);
	\draw[dashed] (m.center) -- (x1.center);
	\draw[dashed] (m.center) -- (x2.center);
\end{tikzpicture}
\caption{a.~A digraph, and b.~the corresponding bounds on resources.}
\label{fig:setsOfBounds}
\end{figure}

We give now a simple illustration of why sets of bounds enable to discard more paths than single bounds. 
Consider the example on Figure~\ref{fig:setsOfBounds}.a, where we have an $o$-$w$ path $Q$, and three $w$-$d$ paths $R$, $R'$, and $R''$. 
Resources, which belong to $\R^2$ endowed with the componentwise sum and order, are indicated on Figure~\ref{fig:setsOfBounds}.a. 
The resources of $w$-$d$ paths are indicated by crosses on Figure~\ref{fig:setsOfBounds}.b.
Consider a situation where $\rmeas\big((\re^1,\re^2)\big) $ is equal to $1$ if and only if $\max(q^1,q^2) > 2$. 
There is no feasible $o$-$d$ path starting by $Q$. 

Suppose first that we are using single bounds as in the usual approach, i.e., $B_w$ contains a unique element $b_w$. Recall that $b_w$ is then such that $b_w\preceq q_S$ for every $w$-$d$ path $S$. In such a case, $b_w \preceq \re_R \meet \re_{R'} \meet \re_{R''} = (1,1)$. Hence, $\rmeas(\re_Q\rplus b_w) = 0$ and the path $P$ is not discarded at Step~\ref{step:LBtest} of Algorithm~\ref{alg:enumeration}. 
Suppose instead that we use the two bounds $b_1 = \re_R \meet \re_{R'}$ and $b_2 = \re_{R''}$, which are indicated by circles on Figure~\ref{fig:setsOfBounds}.b. We have then $\rmeas(\re_Q \rplus b_1) = \rmeas(\re_Q \rplus b_2) = 1$, and the path $Q$ is discarded at Step~\ref{step:LBtest}.
This example is very simple, but it is the same mechanism that is in work in the general case.


The first author~\citep{parmentier2015algorithms} introduced a procedure which, given a size $\kappa$ in input, builds lower bounds sets $B_w$ of size $\kappa$.
Larger sets of bounds $B_w$ enable to get larger lower bounds, and hence to discard more paths. 
However, larger sets of bounds also mean a longer preprocessing is required to compute the bounds. 
Hence, the parameter $\kappa$ is chosen to obtain a tradeoff between the quality of the bounds and the time needed to compute them.
Regarding the practical choice of $\kappa$ for the crew pairing pricing subproblem, the following rule of thumbs ensures good results in practice: Use $\kappa = 1$ if there are fewer than $100$ vertices, $\kappa = 50$ if there are fewer than $300$ vertices, $\kappa = 150$ if there are fewer than $1,500$ vertices, and $\kappa = 250$ if there are more.

When solving \eqref{eq:MasterProblem}, the preprocessing, which actually consists in building an ``extended graph'', is done once and for all: The same extended graph is used each time we solve the pricing subproblem. We can thus work with larger sets of bounds than independent resolutions of the \MRCSP would have allowed.


\section{Integrated problem}
\label{sec:integrated_problem}

\subsection{Problem formulation}
If a crew changes airplane during a connection between two flight legs $\ell$ and $\ell'$, its members need time to cross the airport between the arrival of $\ell$ and the departure of $\ell'$. This is not possible if the time between the arrival of $\ell$ and the departure of $\ell'$ is too short. A \emph{short connection} is an ordered pair $(\ell,\ell')$ of flight legs that can be operated by a crew only if $\ell$ and $\ell'$ are operated by the same airplane. Due to short connections, aircraft routing and crew pairing are linked. Given the collection of all short connections, the integrated problem consists in finding a solution of the aircraft routing problem of Section~\ref{sec:ar} and a solution of the crew pairing problem of Section~\ref{sec:crewPairing} such that whenever a short connection is used in a pairing, it is also used in a route of an airplane.

\subsection{Modeling as an integer program}

The solutions $\ve  x$ of~\eqref{eq:ARprogram} and $\ve y$ of~\eqref{eq:MasterProblem} provide a solution of the integrated problem if and only if
\begin{equation}\label{eq:shortConnectionConstraint}
\sum_{\pairing \in\pairingSet_{\alpha}}y_{\pairing} \leq \sum_{a \in A_{\alpha}} x_{a}
\end{equation} for every short connection $\alpha=(\ell,\ell')$,
where we denote by $A_{\alpha}$ (resp. $\pairingSet_{\alpha}$) the set of arcs (resp. pairings) using the short connection $\alpha$. For any feasible solution of the aircraft routing problem, there is a solution of the crew pairing problem compatible with it since there is no constraint on the number of crews, but solving the two problems simultaneously allows to spare these additional crews and to reduce the costs, as explained in Section~\ref{subsec:review}. The integrated problem aims at performing this task and is thus modeled by the following integer program

\begin{equation}\tag{Int}\label{eq:IntegratedProblem}
	\begin{array}{rll}
	\min & \displaystyle \sum_{\pairing \in \pairingSet} c_{\pairing}y_{\pairing} \\
	\mathrm{s.t.}& \mbox{$\ve  x$ satisfies constraints of~\eqref{eq:ARprogram}}\\
& \mbox{$\ve y$ satisfies constraints of~\eqref{eq:MasterProblem}}\\
& \mbox{$\ve x$ and $\ve y$ satisfy constraints~\eqref{eq:shortConnectionConstraint} for all short connections $\alpha$.}
\end{array}
\end{equation}

\subsection{A cut generating approach}

As we will see in the numerical results, the aircraft routing and crew pairing solution schemes introduced solve most of our industrial instances to optimality in a few hours.
It is therefore natural to test the ability of a simple combination of these approaches to tackle with the integrated problem. Instead of solving directly Program~\eqref{eq:IntegratedProblem}, we adopt a cut generating approach using the methods proposed in the previous sections in a rather independent way.


Let $\scSet(\ve y)$ denote the set of short connections used in a solution $\ve y$ of \eqref{eq:MasterProblem}. Given a feasible solution $\ve y$ of~\eqref{eq:MasterProblem}, if there is no feasible solution $\ve x$ of~\eqref{eq:ARprogram} satisfying \eqref{eq:shortConnectionConstraint}, then any solution $\ve y'$ such that $\scSet(\ve y) \subseteq \scSet(\ve y')$ leads to a more constrained \eqref{eq:ARprogram}, and hence to a similar infeasibility. To avoid such solutions in \eqref{eq:MasterProblem}, we set $S=\scSet(\ve y)$ and add the constraint
\begin{equation}\label{eq:shortConnectionCuts}
	\displaystyle \sum_{\pairing \in \pairingSet } |p \cap \scSet| y_{\pairing} \leq |\scSet|-1,
\end{equation}
where $|p\cap S|$ denotes the cardinality of $\{\alpha\in S\colon p\in\pairingSet_{\alpha}\}$. It prevents a solution to use all short connections in $S$ but does not restrict otherwise the set of solutions.

We can now describe the algorithm for the integrated problem. The algorithm maintains a set $\scConSet$ of short connection cuts. Initially, $\scConSet$ is empty. The following steps are repeated.
\begin{enumerate}[label=(\roman*)]
	\item\label{step:CP} Solve \eqref{eq:MasterProblem} with additional constraints \eqref{eq:shortConnectionCuts} for $\scSet\in\scConSet$. Let $\ve y^*$ be the optimal solution. 
	\item Solve \eqref{eq:ARprogram} with the additional constraints \eqref{eq:shortConnectionConstraint}. \label{step:AR}
	\begin{itemize}
	 	\item If it is feasible, then stop (we have found the optimal solution of \eqref{eq:IntegratedProblem}).
	\item Otherwise, add $\scSet(\ve y^*)$ to $\scConSet$ and go back to \ref{step:CP}.
	\end{itemize} 
\end{enumerate}
Because of the cuts added along the algorithm, a solution $\ve y^*$ is considered at most once. The number of solutions to the crew pairing problem being finite, the cut generation algorithm terminates after a finite number of iterations. The solutions of the last call to \ref{step:CP} and  \ref{step:AR} form an optimal solution to \eqref{eq:IntegratedProblem}: at each iteration, the only solutions to \eqref{eq:MasterProblem} that are forbidden by the additional constraints \eqref{eq:shortConnectionCuts} are not feasible for \eqref{eq:IntegratedProblem} and at the last iteration, $\ve y^*$ is the optimal solution of a relaxation of \eqref{eq:IntegratedProblem}.

In practice, the algorithm does not converge after thousands of iterations on industrial instances. We therefore replace $|S|-1$ by $\gamma |S|$ with $\gamma < 1$ in the constraints \eqref{eq:shortConnectionCuts}, losing the optimality of the solution returned. 
During the first iteration, there is no additional constraint \eqref{eq:shortConnectionCuts}: the crew pairing problem \eqref{eq:MasterProblem} is therefore not constrained by the aircraft routing.
It is therefore a relaxation of the integrated problem, and its optimal solution provides a lower bound on the optimal solution of the integrated problem. We use this lower bound to evaluate the quality of the solution of the integrated problem returned by the algorithm.
Numerical experiments in Section \ref{sec:experimentalResults} show that $\gamma=0.9$ is a good compromise: We obtain near optimal solutions after a few dozens of iterations.

\section{Experimental results}
\label{sec:experimentalResults}


\subsection{Instances} 
\label{sub:instances}

Table~\ref{tab:subfleetInstance} describes six industrial instances of Air France. Each instance contains the legs of a fleet on a weekly horizon. The first two columns provide the name of the instance and the number of legs it contains. Columns ``Airplane connect.''~and ``Airplanes'' respectively give the number of connections that can be done by airplanes, i.e.,~the number of ordered pairs of legs $(\ell,\ell')$ that can be operated consecutively in a route, and the number $n^{\mathrm{a}}$ of airplanes available. Columns ``Crew connect.''~and ``Crew pairings'' respectively provide the number of connections that can be taken by crews, and the order of magnitude of the number of pairings in a good solution. Finally, column ``Short connect.''~gives the number of short connections available. These instances are large: For instance, the largest instance considered by 
\citet{mercier2005computational} has 707 legs,
and the largest instance for the integrated problem in the literature~\citep{weide2010iterative} has 750 legs.  

The instance A318-9 (resp.~A320-fam) contains the legs of the A318 and A319 (resp.~A318, A319, A320, and A321) instances, as well as a few extra ``fictitious'' legs. As Air France's crews can operate legs on planes of different subfleets on the same pairing, there is a common crew pairing problem for each of the instance A318-9 and~A320-fam. 
On the contrary, the subfleet of an airplane is fixed: An A318 does not become an A319. Hence, solving aircraft routing for multiple subfleets together consists in solving one separate problem for each subfleet. This is what we do when we solve aircraft routing instances within the integrated problem solution scheme for instances A318-9 and~A320-fam.


\begin{table}
	\begin{center}
	\begin{tabular}{|lrrrrrr|}
		\hline
		Instance & \multicolumn{1}{c}{Legs} & Airplane & Airplanes & \multicolumn{1}{c}{Crew} & \multicolumn{1}{c}{Crew} & \multicolumn{1}{c|}{Short} \\
		& & connect. & & connect. & \multicolumn{1}{c}{pairings ($\simeq$)} & connect. \\
		\hline
		A318 & 669 & 39,564 & 18 & 3,742 & 130 & 1,230 \\
		A319 & 957 & 45,901 & 41 & 3,738 & 240 & 996 \\
		A320 & 918 & 49,647 & 45 & 3,813 & 280 & 1,103 \\
		A321 & 778 & 29,841 & 25 & 3,918 & 165 & 1,006 \\
		A318-9 & 1,766 & -- & (59) & 8,070 & 350 & 2,226\\
		A320-fam & 3,398 & -- & (129) & 21,563 & 690 & 4,398 \\
		\hline
	\end{tabular}
	\end{center}
	\caption{Air France industrial instances}
	\label{tab:subfleetInstance}
\end{table}

\subsection{Experimental setting} 
\label{sub:experimental_setting}

All the numerical experiments are performed on a server with 128 GB of RAM and 12 cores at 2.4 GHz. \texttt{CPLEX 12.1.0} is used to solve all linear and integer programs. The algorithms are not parallelized.

\subsection{Aircraft Routing} 
\label{sub:aircraft_routing}

\Cref{tab:ARres} provides the results for aircraft routing. 
The first column gives the name of the instance.
The next two ones give results for the aircraft routing problem on its own.
Column ``Uncons.~CPU time (mm:ss)'' gives the time needed to solve \eqref{eq:ARprogram}. 
Column ``Optim.~CPU time (mm:ss)'' gives the time needed to find an optimal solution of the optimization problem that consists in finding the minimum number of airplanes needed to operate the instance: We use the left-hand side of \eqref{eq:number_airplanes} as objective. 
Note that this problem has not been mentioned previously in the paper.
The solution scheme for the integrated problem in \Cref{sec:integrated_problem} solves \eqref{eq:ARprogram} with additional constraints~\eqref{eq:shortConnectionConstraint}. 
The last two columns provide numerical results for this constrained version.
 On all but the last iterations of the integrated problem scheme, aircraft routing is infeasible. 
Column ``Infeas.~CPU time'' provides the time needed to solve the penultimate iteration, which is infeasible, and column ``Feas.~CPU time '' the last iteration, which is feasible.
 The typical computing time is a few dozens of seconds on industrial instances. The longest constrained feasible version requires a few minutes.
The optimization version in the second column is typically one order of magnitude faster than the one obtained by \citet[Tables 10 and 11]{khaled2018compact} on instances with similar numbers of legs. 
However, this last statement must be taken with care as the structure of the instances (number of airplanes, number of days of planning, etc.) is very different. 
This improved performance is likely due to the fact that our relaxation has a better relaxation than their one, as we prove in Proposition~\ref{prop:compareKhaled} in Appendix~\ref{sec:proofs} and that our formulation has less symmetry.



\begin{table}
	\begin{center}
	\begin{tabular}{|l|rr|rr|}
		\hline
		
		& \multicolumn{2}{|c|}{\centering \eqref{eq:ARprogram} alone}
		& \multicolumn{2}{c|}{\centering \eqref{eq:ARprogram} within \eqref{eq:IntegratedProblem}} \\
		Instance  
		& \multicolumn{1}{|p{2.5cm}}{\centering Uncons.~CPU time (mm:ss)}
		& \multicolumn{1}{p{2.5cm}|}{\centering Optim.~CPU time (mm:ss)}
		& \multicolumn{1}{p{2.5cm}}{\centering Infeas.~CPU time (mm:ss)}
		& \multicolumn{1}{p{2.5cm}|}{\centering Feas.~CPU time (mm:ss)}
		\\
		\hline
		A318 & 00:17 & 00:58  & 00:14 & 01:35 \\
		A319 & 00:16 & 01:05  & 00:22 & 00:19 \\
		A320 & 01:02 & 03:55  & 00:35 & 13:28 \\
		A321 & 00:16 & 01:03  & 00:23 & 00:19 \\
		\hline
	\end{tabular}
	\end{center}
	\caption{Aircraft routing results}
	\label{tab:ARres}
\end{table}

\subsection{Crew Pairing} 
\label{sub:crew_pairing}


\Cref{tab:CPresults} provides the results for crew pairing. 
All instances are solved to optimality.
The first column of \Cref{tab:CPresults} gives the name of the instance. The next column provides the value of $\kappa$ determined using the rule of thumb of Section \ref{sub:enhanced_bounds_computations_for_the_column_generation_context} and needed by the algorithm building the sets $B_v$. 
Column ``Col.~Gen.~Iter'' provides the number of iterations in the column generation, and column ``Pricing time'' the percentage of time spent in the pricing subproblem. 
This pricing time includes the time needed by the computation of the sets $B_v$ and by the enumeration algorithm (Algorithm~\ref{alg:enumeration}). Columns ``LP time'' and ``MIP time'' indicate the percentage of the total CPU time spent in Algorithm~\ref{alg:CG} solving Step~\ref{step:master}, and solving Steps~\ref{step:firstUB} and~\ref{step:final}.
The last column gives the total time needed by the algorithm.
On all these instances, the integrality gap does not exceed 0.01\%. This explains the fast resolution of Step~\ref{step:final}.

\begin{table}
	\begin{tabular}{|l|r|rrrr|r|}
	\hline
		Instance 
		& \multicolumn{1}{p{1cm}|}{\centering $\kappa$} 
		& \multicolumn{1}{p{1.4cm}}{\centering Col.~Gen. Iter} 
		& \multicolumn{1}{p{1cm}}{\centering Pricing time} 
		& \multicolumn{1}{p{1cm}}{\centering LP time} 
		& \multicolumn{1}{p{1cm}|}{\centering MIP time} 
		& \multicolumn{1}{p{2cm}|}{\centering Total~time (hh:mm:ss)}\\
	\hline
A318    & 150 & 394  & 86.60\% & 13.34\% & 0.05\% & 01:21:22 \\ 
A319    & 150 & 264  & 60.66\% & 39.14\% & 0.15\% & 00:10:47 \\ 
A320    & 150 & 226  & 74.54\% & 25.20\% & 0.20\% & 00:08:35 \\ 
A321    & 150 & 382  & 65.82\% & 32.60\% & 1.25\% & 00:33:51 \\ 
A318-9  & 150 & 867  & 69.71\% & 30.21\% & 0.07\% & 05:43:00 \\ 
A320fam & 250 & 2,166 & 43.28\% & 56.62\% & 0.10\% & 104:05:59 \\ 
	\hline
	\end{tabular}
	\caption{Crew pairing results  -- Instances are solved to optimality}
	\label{tab:CPresults}
\end{table}

	\begin{rem}
		One may be tempted to stop the column generation before convergence to exchange quality for speed.
		Unfortunately, and this is a limit of our method, if we stop the column generation before convergence, the solution found by the MIP solver at Step~\ref{step:firstUB} is poor, and Step~\ref{step:completeWithPathGap} is not tractable in practice.
	\end{rem}


\subsection{Integrated problem} 
\label{sub:integrated_problem}



\Cref{tab:IPresults} provides the results for the integrated problem. The constraint strength parameter $\gamma$ of the end of Section~\ref{sec:integrated_problem} is equal to 0.9, and the bounds sets size $\kappa$ is equal to $150$. 
The first column provides the instance solved.
Columns ``Integ.~steps'' provides the number of steps of the integrated problem algorithm of Section~\ref{sec:integrated_problem} before convergence. 
Column ``CG it.~total'' provides the total number of column generation iterations realized on the successive integrated problem algorithms steps. Column ``\eqref{eq:MasterProblem} CG time'' provides the proportion of the total CPU time spent in the column generation, i.e., solving the pricing subproblem and the linear relaxation of the master problem,
and column ``\eqref{eq:MasterProblem} MIP time'' the proportion spent solving the integer version of the crew pairing master problem.
The column ``\eqref{eq:ARprogram} time'' provides the proportion spent solving aircraft routing integer program \eqref{eq:ARprogram}. The column ``Sho.~Con.'' gives the number of short connections in the final solution. The linear relaxation of the crew pairing master problem \eqref{eq:MasterProblem} with no short connection constraint is used as the lower bound on the cost of an optimal solution. The gap provided is between the cost of the solution returned and this lower bound. Finally, the last column provides the total CPU time needed by the algorithm. Only instance A320-fam could not be solved, as the algorithm had not converged after one week of computing time.

\begin{table}
\begin{outdent}
	\begin{tabular}{|l|rrrr|r|rr|r|}
	\hline 
		Instance 
		& \multicolumn{1}{p{1cm}}{\centering Integ. steps} 
		&  \multicolumn{1}{p{1.1cm}}{\centering CG~it. total} 
		& \multicolumn{1}{p{1.35cm}}{\centering \eqref{eq:MasterProblem}~CG time} 
		& \multicolumn{1}{p{1.60cm}|}{\centering \eqref{eq:MasterProblem}~MIP time} 
		& \multicolumn{1}{p{1cm}|}{\centering \eqref{eq:ARprogram} time} 
		& \multicolumn{1}{p{1cm}}{\centering Sho. Con.} 
		& \multicolumn{1}{c|}{Gap} 
		& \multicolumn{1}{p{2.3cm}|}{\centering Total~time (hh:mm:ss)} \\ 
		\hline
		A318 & 6 &  460 & 95.53\% & 2.56\% & 1.91\% & 323 & 0.0002\% & 01:53:47\\
 		A319 & 4 &  343 & 76.99\% & 13.27\% &  9.74\% & 448 & 0.0013\% & 00:20:18 \\
		A320 & 2 &  240 & 24.24\% & 39.38\% & 36.38\% & 436 & 0.0017\% & 00:38:36 \\
		A321 & 2 &  380 & 96.60\% &  2.53\% &  0.88\% & 413 & 0.0074\% & 00:29:18 \\
		A318-9 & 2 &  915 & 97.66\% & 1.71\% & 0.63\% & 790 & 0.0008\% & 06:34:31 \\
		A320-fam& \multicolumn{8}{|c|}{Stopped after one week} \\
		\hline
	\end{tabular}
		\caption{Numerical results on integrated problem}
		\label{tab:IPresults}
\end{outdent}
\end{table}

We emphasize the fact that the solution returned by the approximate algorithm is almost optimal. Practically speaking, the gap obtained is less than or equal to 0.01\%. 
 The computation time needed to obtain a near optimal solution of the integrated problem is of the same order of magnitude than the time needed to obtain a solution of the crew pairing problem in \Cref{tab:CPresults}. Solving aircraft routing and crew pairing sequentially strongly constrains the solution: Indeed, when solved in an integrated fashion, around half of the connections in the solution are short connections.

 \begin{table}
 \begin{outdent}
 \begin{tabular}{|l|rrrrr|}
 \hline
 Instance & 318 & 319 & 320 & 321 & 318-9 \\
 \hline
 Cost~reduction & 0.08\% & 0.16\% & 0.33\% & 0.31\% & 0.31\% \\  
 CPU time ratio & 10.6$\times$ & 3.0$\times$ & 5.6$\times$ & 2.0$\times$ & 5.1$\times$ \\ 
 \hline
 \end{tabular}
 \end{outdent}
 \caption{Sequential versus integrated resolution of aircraft routing and crew pairing}
 \label{tab:seqVsInt}
 \end{table}
 Finally, Table~\ref{tab:seqVsInt} compares the sequential approach, where aircraft routing is solved first, and then crew pairing, to the integrated approach of Section~\ref{sec:integrated_problem}. Line ``Cost~reduction'' provides the percentage by which the cost is reduced when using the integrated approach, and line ``CPU time ratio'' the increase in computing time. On average, using the integrated approach enables to reduce the costs by 0.25\%, and computing time is 5.1 times longer. This reduction of cost on our instances is smaller than what is mentioned in the literature: $5\%$ on average according to \citet{cordeau2001benders}, and $1.6\%$ according to \citet{papadakos2009integrated}. We believe that this comes from the fact that our instances are larger: Adding new connections have a stronger impact when few connections are available. The increase in computing time is mainly due to the fact that the crew pairing is longer to solve in the integrated approach due to the addition of short connections.

\subsection{Industrial relevance}
\label{sub:discussion}
To be usable in an industrial context, the computing time of the solvers must not exceed eight hours, which represent one night of computing time. Our algorithms enable to solve to near optimality instances of the integrated problem with up to $1,766$ legs within this time constraint. 
Our solution scheme therefore enables to deal practically with instances larger than those in the literature -- the largest instances in the literature \citep{weide2010iterative} have 750 legs. 
These performances have been made possible by our pricing subproblem algorithm. Further improvements to deal with larger instances cannot be done by only working on the pricing subproblem. Indeed, we can see in Table~\ref{tab:CPresults} that most of the CPU time on instance A320-fam is spent in the simplex algorithm.

\subsection{Focus on the pricing subproblem} 
\label{sub:focus_on_the_pricing_subproblem_choice_of_}

\begin{table}
	\begin{center}
	\begin{tabular}{|lrrr|}
		\hline
		Instance 
		& \multicolumn{1}{c}{Legs} 
		& \multicolumn{1}{c}{Connections} 
		& \multicolumn{1}{c|}{Pairings ($\simeq$)} \\
		
		\hline
		CP50 & 290 & 1,006 & 50\\
		CP70 & 408 & 1,705 & 70\\
		CP90 & 516 & 2,490 & 90\\
		\hline
	\end{tabular}
	\end{center}
	\caption{Medium size artificial crew pairing instances}
	\label{tab:CPextractedInstances}
\end{table}


The key element in the performance of our approach is the performance of Algorithm~\ref{alg:enumeration}. 
We compare it in this section to the algorithms previously used.
As our industrial instances are too large to be solved using these algorithms, we introduce in Table~\ref{tab:CPextractedInstances} smaller instances, which we have built by considering only a subset of the legs of the instance A318. 
Columns ``Legs'' and ``Connections'' respectively provide the number of legs and connections in the instances, and column ``Pairings'' the approximate number of pairings in a solution. \Cref{tab:CPsubproblemAlgorithmsResults} provides results on the performance of the pricing subproblem scheme on these instances. 
Its first column gives the instance solved.
The next one provides the algorithm used. 
Parameter $\kappa$, introduced in Section~\ref{sub:enhanced_bounds_computations_for_the_column_generation_context}, gives the size of the lower bounds sets for algorithms using bounds.
The next three columns provide statistics on the resource constrained shortest path (RCSP) algorithms. 
As mentioned in Section~\ref{sub:modeling_the_pricing_subproblem}, seven RCSP instances are solved for each pricing subproblem, one for each period of four consecutive days.
The statistics are averaged on all the instances solved along the column generation.
Column ``RCSP~iter av.~nb'' provides the average number of iterations of the RCSP algorithm, ``Cut~Dom.''~provides the proportion of paths cut at Step~\ref{step:DomTest}, the remaining being cut at Step~\ref{step:LBtest}. 
Column ``RCSP time'' provides the average time needed to solve one RCSP instance.
Column ``Pricing subproblem'' provides the proportion of the total computing time spent solving the pricing subproblem, and the last column gives the total computing time of the crew pairing solution scheme.

We can see that the use of bounds enables a huge speed-up with respect to the usual algorithms.  
This speed-up is required to deal with instances with more than 500 legs.
Two elements explain this speed-up.
First, the condition at Step~\ref{step:LBtest} enables to discard many paths:
In Algorithm~\ref{alg:enumeration}, more than 90\% of the paths discarded are discarded at Step~\ref{step:LBtest} and not at Step~\ref{step:DomTest}.
Second, running Algorithm~\ref{alg:enumeration} with $c(q_P)$ instead of $\min \{c(\re_P \rplus b)\colon b \in B_v, \rmeas(\re_P \rplus b)=0\}$ makes it much slower. Hence, $\min \{c(\re_P \rplus b)\colon b \in B_v, \rmeas(\re_P \rplus b)=0\}$  seems numerically to be a better evaluation of how $P$ is promising than $c(q_P)$, and enables Algorithm~\ref{alg:enumeration} to find good solutions faster than the label correcting algorithm, and also faster than the label setting algorithm that does not use keys at all for choosing the next path to consider.

\begin{table}
	\begin{outdent}
\begin{tabular}{|l|cr|rrr|rr|}
	\hline
		Instance 
		& Algorithm\ 
		& \multicolumn{1}{p{1cm}|}{\centering $\kappa$} 
		& \multicolumn{1}{p{1.7cm}}{\centering RCSP~iter av.~nb.} 
		& \multicolumn{1}{p{1cm}}{\centering Cut Dom.} 
		& \multicolumn{1}{p{2.2cm}|}{\centering RCSP~time av~(mm:ss.ff)} 
		& \multicolumn{1}{p{1cm}}{\centering Pricing time} 
		& \multicolumn{1}{p{1.9cm}|}{Total~time (hh:mm:ss)} \\
	\hline

CP50 & Label setting & -- & 1.020e+04 &  --  & 00:00.56 & 97.55\% & 00:04:38\\ 
CP50 & Label correcting   & -- & 1.308e+04 &  --  & 00:01.28 & 97.38\% & 00:11:37\\ 
CP50 & Algorithm~\ref{alg:enumeration} & 1   & 2.326e+03 & 6.89\% & 00:00.03 & 75.28\% & 00:00:23\\ 
CP50 & Algorithm~\ref{alg:enumeration} & 10  & 4.914e+02 & 4.01\% & 00:00.02 & 59.87\% & 00:00:17\\ 
CP50 & Algorithm~\ref{alg:enumeration} & 100 & 2.033e+02 & 5.03\% & 00:00.04 & 77.06\% & 00:00:33\\ 
\hline 
CP70 & Label setting & -- & 5.644e+04 &  --  & 00:11.49 & 99.52\% & 05:07:05\\ 
CP70 & Label correcting   & -- & 7.730e+04 &  --  & 00:17.16 & 99.56\% & 07:28:22\\ 
CP70 & Algorithm~\ref{alg:enumeration} & 1   & 9.208e+03 & 7.69\% & 00:00.24 & 90.61\% & 00:04:41\\ 
CP70 & Algorithm~\ref{alg:enumeration} & 10  & 1.994e+03 & 4.28\% & 00:00.04 & 58.48\% & 00:01:12\\ 
CP70 & Algorithm~\ref{alg:enumeration} & 100 & 8.007e+02 & 5.77\% & 00:00.07 & 77.43\% & 00:01:43\\ 
\hline
CP90 & Label setting & -- & 9.779e+04 &  --  & 00:40.71 & \multicolumn{2}{c|}{Stopped after 48h}\\ 
CP90 & Label correcting   & -- & 2.007e+05 &  --  & 01:42.87 & \multicolumn{2}{c|}{Stopped after 48h}\\ 
CP90 & Algorithm~\ref{alg:enumeration} & 1   & 5.000e+04 & 9.81\% & 00:05.98 & 98.86\% & 02:56:33\\ 
CP90 & Algorithm~\ref{alg:enumeration} & 10  & 9.966e+03 & 5.88\% & 00:00.34 & 81.86\% & 00:12:36\\ 
CP90 & Algorithm~\ref{alg:enumeration} & 100 & 4.377e+03 & 5.60\% & 00:00.25 & 77.98\% & 00:10:28\\ 
\hline
A318 & Label setting & -- & 1.319e+05 &  --  & 00:53.01 & \multicolumn{2}{c|}{Stopped after 48h}\\ 
A318 & Label correcting   & -- & 3.802e+05 &  --  & 01:36.04 & \multicolumn{2}{c|}{Stopped after 48h}\\ 
A318 & Algorithm~\ref{alg:enumeration} & 1   & 7.161e+04 & 8.99\% & 00:08.61 & 97.87\% & 05:35:42\\ 
A318 & Algorithm~\ref{alg:enumeration} & 10  & 5.472e+04 & 6.62\% & 00:05.97 & 96.02\% & 05:06:47\\ 
A318 & Algorithm~\ref{alg:enumeration} & 100 & 2.549e+04 & 3.72\% & 00:01.65 & 86.97\% & 01:32:50\\ 

	\hline
	\end{tabular}

	\end{outdent}
	\caption{Relative performance of pricing subproblems algorithms}
	\label{tab:CPsubproblemAlgorithmsResults}
\end{table}

\section{Conclusion} 
\label{sec:conclusion}
We have proposed a compact integer program for aircraft routing. Its main strength is its ease of implementation. Numerical results show that it can deal with large industrial instances in at most a few minutes, even when optimization versions are considered. We have used a resource constrained shortest path algorithm recently introduced by the first author for the crew pairing column generation pricing subproblem. This algorithm leverages on the lattice ordered monoid structure of the resource set to build efficient lower bounds. Practically, this enables to solve to optimality very large industrial crew pairing instances. As a side result, we have combined these aircraft routing and crew pairing solution schemes in a cutting plane approach to the integrated problem. The resulting algorithm solves to near optimality large industrial instances of the integrated problem.

\section*{Acknowledgments} 
\label{sec:acknowledgments}
We thank the reviewers for their thorough reading and all their useful comments and suggestions that helped us improve the paper.
We are grateful to Air France which partially supported the project. We are especially thankful to Alexandre Boissy, who initiated the project, and to Mathieu Sanchez and Mohand Ait Alamara who helped us with the implementation in Air France softwares. We also thank Sourour Elloumi for pointing out the correct reference for Lemma~\ref{lem:baldacci}.


\bibliographystyle{plainnat}
\bibliography{air}

\appendix


\section{Proofs}
\label{sec:proofs}

\newcommand{\redcost}{\tilde{c}}


\begin{proof}[Proof of Proposition~\ref{prop:AR}]
The fact that a feasible solution of the aircraft routing problem induces a collection $\C$ as in the statement is obvious. Let us prove the other direction, which is almost as easy.

Consider a collection $\C$ as in the statement. Each cycle provides a route, possibly of several weeks. We show now that the solution consisting of these routes is feasible. By construction of the graph, the maintenance requirement is satisfied. Moreover, the number of times a cycle intersects $A_0$ is an upper bound on the number of airplanes required to operate the corresponding route: The number of times it intersects $A_0$ is equal to the number of weeks this cycle lasts, as illustrated on Figure~\ref{fig:timeline}. Thus the number of arcs selected in $A_0$ by the whole collection is an upper bound on the number of airplanes required to operate the solution. Since this number is at most $n^{\mathrm{a}}$, the solution is feasible.
\end{proof}

\begin{proof}[Proof of Lemma~\ref{lem:baldacci}]
Consider the following integer program in the standard form:
\begin{equation}\label{eq:IPgap}
	\begin{array}{rl}
	\min \enskip & \cc^T\zz \\
	\mathrm{s.t.} \enskip & A\zz =\bb \\
	& \zz \in \Z_+.
	\end{array}
\end{equation}
By the theory of the simplex algorithm, the program~\eqref{eq:IPgap} can be written in the form
$$
\begin{array}{rl}
	\min \enskip & \ds{\bar v+\sum_{i\in N}\tilde c_iz_i} \\
	\mathrm{s.t.} \enskip & A\zz =\bb \\
	& \zz \in \Z_+,
	\end{array}
$$
where the $\tilde c_i$ -- the reduced costs -- are all non-negative, and where $N$ is the complement of the considered basis (all indices $i$ such that $\tilde c_i>0$ are in $N$). Consider now an optimal solution $\zz^*$ of the program~\eqref{eq:IPgap}.
We have $$\bar v+\sum_{i\in N}\tilde c_iz_i^*=\cc\cdot\zz^*\leq \UB .$$ The non-negativity of the $\tilde c_i$'s and the $z_i^*$'s implies the following inequality for every $i\in N$:
$$z_i^*\leq \frac{\UB-\bar v}{\tilde c_i}.$$ In particular, if $\tilde c_i>\UB-\bar v$, we necessarily have $z_i^*<1$, i.e., $z_i^*=0$.
\end{proof}

\begin{proof}[Proof of \Cref{prop:MRCSPalgoConvergence}]
We first prove that the algorithm terminates after a finite number of iterations. By induction on the iterations, we see that there is never two distinct elements in $\mathsf{L}$ such that one is a subpath of the other: There is no such two elements when the algorithm starts, and if there is no such two elements at a given iteration, there cannot be such two elements at the next iteration either. In particular, when a path $P$ leaves $\mathsf{L}$ at some iteration, it prevents the presence of a subpath of it in $\mathsf{L}$ at the current iteration. Thus, because of the update rule, $P$ cannot be added again to $\mathsf{L}$ in a subsequent iteration. It implies that a given path $P$ is considered at most once by the algorithm and, as there is a finite number of paths in an acyclic digraph, we get the sought conclusion regarding the termination of the algorithm.


We prove now the part of the statement regarding the cost of $c_{od}^{UB}$ at the end of the algorithm.
At any step of the algorithm, $c_{od}^{UB}$ is either equal to $+\infty$ or to the cost of an $o$-$d$ path $P$ such that $\rmeas(\re_{P}) = 0$. Therefore, if there is no feasible solution, then $c_{od}^{UB}$ is never updated, and equal to $+\infty$ at the end of the algorithm. Suppose now that there is a feasible solution, and let $P$ be a feasible $o$-$d$ path of minimum cost. By definition of $P$ and due to the update rule, we have $\rcost(\re_{P}) \leq c_{od}^{UB}$ at the end of the algorithm. Suppose for a contradiction that this inequality is strict. Given the update rule of $c_{od}^{UB}$, this means that neither $P$ nor a feasible $o$-$d$ path $Q$ dominating $P$ has been considered present in $\mathsf{L}$. Let $P'$ be the longest subpath of $P$, with origin $o$, such that $P'$, or a path dominating $P'$ with the same origin and destination as $P'$, has been present in $\mathsf{L}$. Denote by $v$ the destination of $P$.
Among all $o$-$v$ paths present in $\mathsf{L}$ at some time during the algorithm and that dominate $P'$, pick a path $Q'$ that is non-dominated by the others.

The test of Step~\ref{step:DomTest} is necessarily satisfied by $Q'$ because, by definition of $Q'$, there is no $o$-$v$ path in $\mathsf{L}_{v}^{\mathrm{nd}}$ dominating $Q'$ when $Q'$ is considered by the algorithm. We now prove that the test of  Step~\ref{step:LBtest}  is also necessarily satisfied. Indeed, let $P''$ be the $v$-$d$ subpath such that $P=P'+P''$, and $b$ a bound in $B_{v}$ such that $b \rleq \re_{P''}$. (Here, $P'+P''$ means that we append $P''$ to $P'$.) We have $\re_{Q'} \rplus b \rleq \re_{P'} \rplus b \rleq \re_{P'} \rplus \re_{P''} = \re_{P}$. We therefore have, when $Q'$ is considered, $\rmeas(\re_{Q'} \rplus b) \leq \rmeas(\re_{P}) = 0$, and $\rcost(\re_{Q'} \rplus b) \leq \rcost(\re_{P}) < c_{od}^{UB}$, where the last inequality relies on the fact that $c_{od}^{UB}$ is non-increasing along the algorithm. Hence, the test of Step~\ref{step:LBtest} is satisfied. Therefore, $Q'$ satisfies the tests of  Steps~\ref{step:LBtest} and~\ref{step:DomTest}, and is added to $\mathsf{L}$ whatever the combination of these two tests the algorithm uses. Since $Q'$ is non-dominated by other $o$-$v$ paths considered by the algorithm, it is extended in a subsequent iteration. Denote by $a$ the arc in $P$ that follows $P'$ and denote by $w$ the head of $a$. Since $q_{Q'}\rplus q_a\preceq q_P'\rplus q_a$, either $Q'+a$, or another $o$-$w$ path that dominates $P'+a$ is present once in $\mathsf{L}$, which contradicts the length maximality of $P'$.
\end{proof}

\begin{proof}[Proof of Lemma~\ref{lem:monoid}]
Considering the different cases in the definitions enables to prove that $(0,0)$ is the neutral element of $\rplus$, that $\rplus$ is associative, that $\rleq$ is compatible with $\rplus$, and that $(\rset^{\rmeas}, \rleq)$ is a lattice with meet operator
\begin{align*}
\re \meet \infty &= \re \quad  \\
(n,f) \meet (\tilde{n},\tilde{f}) &= \big(\min(n,\tilde{n}),\min(f,\tilde{f})\big) \\
(n^{b},f^{b},n^{e},f^{e}) \meet (\tilde{n}^{b},\tilde{f}^{b},\tilde{n}^{e},\tilde{f}^{e})  & =
\big(
\min(n^{b},\tilde{n}^{b})
\min(f^{b},\tilde{f}^{b})
\min(n^{e},\tilde{n}^{e})
\min(f^{b},\tilde{f}^{b})
\big)
\end{align*}
and $q \meet \tilde{q} = (0,0)$ for any other combinations.
\end{proof}


\begin{proof}[Proof of Proposition~\ref{prop:pricing}]
Let $p$ be a sequence of flight legs. The definition of $D$ ensures that there is an $o$-$d$ path $P$ whose vertices correspond to the legs in $p$ if and only if $P$ is a pairing that satisfies rules~\ref{rule:length} and~\ref{rule:Paris}. In that case, this path $P$ is unique. Let $p$ be such a pairing and $P$ be the corresponding path.
Let $q_P = \big(q_P^{\rmeas},c_P\big)$ be the resource of $P$.

By definition of the arc resources, $c_P$ is the sum of $c_{(\ell,\ell')} + z_{\ell'}$ for $(\ell,\ell')$ in $P$ and hence in $p$, and we therefore have $c_P = \tilde{c}_p$. Hence $c(q_P) = \tilde{c}_p$.

Let $n^{b}$, $f^{b}$, $n^{e}$, and $f^{e}$ be respectively the number of legs and the flying time of the first duty of $p$, and the number of legs and the flying time of the last duty of $p$. Let $A_P$ be the set of duties of $P$ except the first and the last.
We claim that $q_P^{\rmeas}$ is equal to $\infty$ if and only if there is a duty in $A_P$ that does not satisfy  both rules~\ref{rule:NbLegs} and~\ref{rule:FlyingTime}, and to $(n^{b},f^{b},n^{e},f^{e})$ otherwise.
This result is proved by induction on the number of arcs in $P$. Denoting $a$ the last arc of $P$, and $P'$ the subpath of $P$ obtained by removing~$a$, the induction hypothesis can be applied to $P'$, and the result for $P$ follows by considering the different possible cases for the components in $\rset^\rmeas$ of $\re_P$ and $\re_a$. 
The definition of $\rset^\rmeas$ then ensures that $\rmeas(\re_{P}) = 0$ if and only if rules~\ref{rule:NbLegs} and~\ref{rule:FlyingTime} are satisfied. Hence, $\rmeas(\re_{P}) = 0$ if and only if $p$ is a feasible pairing.
\end{proof}


\section{Aircraft routing MIP}
\label{sec:aircraft_routing_mip}


\citet{khaled2018compact} propose a compact MIP for the tail assignment problem.
Aircraft routing and tail assignment both consist in building the sequences of legs operated by the airplanes of an airline. 
The main difference between them is that in aircraft routing, airplanes are identical and routes do not need to be assigned to airplanes, while in tail assignment airplane specific costs are taken into account and routes are assigned to airplanes.
Aircraft routing and tail assignment problems being fairly similar, MIP formulations for one problem can generally be applied to the other one.

However, the first reason why we do not use Khaled et al.'s formulation on Air France problem is that their formulation does not naturally adapt to this problem.
Indeed,
an important difference between Khaled et al.'s tail assignment and Air France aircraft routing is that, 
in Khaled et al.'s problem, routes are not cyclic, 
while in Air France problem, they are.
By cyclic routes, we mean that sequences of legs are built for a typical week, and that airplanes operate the same sequences of legs week after week in a cyclic way.
And as we have seen in Section~\ref{subsec:ar_formulation}, one cycle can potentially last several weeks.
If Khaled et al.'s compact MIP for tail assignment can easily be adapted to solve a ``non-cyclic'' version of the aircraft routing problem, this is not the case for the ``cyclic'' aircraft routing problem considered at Air France that we introduce in Section~\ref{sec:ar}.
Indeed, adapting it would require introducing many new binary variables to encode how the end of a week cycles with the beginning of the next one.

The second reason is computational.
As cyclic routes are difficult to take into account in their formulation, we compare the formulations on the Air France tail assignment problem.
Both Khaled et al.'s formulation and an adapted version of our formulation \eqref{eq:ARprogram} have been implemented at Air France, and the adapted version of our formulation is now used in practice by Air France to solve its tail assignment.
Our formulation is able to solve the Air France medium haul instances (the tail assignment version of instances A318, A319, A320, and A321 of Table~\ref{tab:subfleetInstance}) to optimality in at most 3 minutes, while the version of Khaled et al. is unable to solve these instances in two hours.
On long-haul instances of Air France, which are easier, their formulation had been previously used and was able to find optimal solutions, but each instance took at least 16 minutes, while our formulation solves each instance in at most 20 seconds.

Actually, the better performances of our formulation are easily explained by theoretical considerations on the linear relaxations. Solvers of MIP are based on branch-and-bound, which crucially relies on the quality of the linear relaxation to discard partial solutions.  
On the Air France tail assignment problem, the linear relaxation of our MIP provides bounds that are non smaller than those provided by the Khaled et al.'s MIP and can be strictly larger, even on very simple and natural examples. 
In the remaining of the appendix, we introduce the Air France tail assignment problem, adapt Khaled et al.'s formulation and \eqref{eq:ARprogram} to that problem, and prove the result mentioned on linear relaxations.

\medskip

We now introduce the \emph{Air France tail assignment problem}.
The input is formed of a given week, a set of airports, a collection $\mathcal{L}$ of flight legs operated between these airports that week, and a set of $\na$ available airplanes. 
Some airports are bases where maintenance can be performed, and each airplane must still spend a night in a maintenance base at least every $\maintenance$ days.
Airplanes are indexed by $j$.
For each airplane $j$ in $[\na]$, let $k_0^j$ be the airport where airplane $j$ is
 at the beginning of the week, 
 and $\delta_0^j$ be the number of days since the last maintenance night of $j$ at the beginning of the week.
For each airplane $j$ and leg $\ell$, 
we have a cost of operating a leg $\ell$ with airplane $j$.
The aim of the tail assignment problem is to build the (non-cyclic) sequence of legs operated by each airplane at minimum cost.

A \emph{tail assignment connection} is therefore a pair $(\ell,\ell')$ of flight legs such that $\ell'$ departs from the arrival airport of $\ell$, and such that the departure time of $\ell'$ minus the departure time of $\ell$ is bounded from below by a given quantity. 
We underline that, if there were airplane connections between a leg $\ell$ at the end of the week and a leg $\ell'$ at the beginning of the week in the aircraft routing problem, there is no such tail assignment connection.
A \emph{tail assignment route $r$ for airplane $j$} is a (non-cyclic) sequence of legs $\ell_1,\ldots,\ell_k$ such that $\ell_i$ departs from $k_0^i$ and any pair of two consecutive legs $(\ell_i,\ell_{i+1})$ is a tail assignment connection.
It satisfies the \emph{maintenance requirement} if, first, supposing that airplane $j$ follows this route, it spends its first night in a maintenance base after at most $\maintenance - \delta_0^j$ days and then spends a night in a maintenance base at least every $\maintenance$ days, and second, if $\ell_k$ arrives in an airport that is not a base, then the last night of the week spent in a base is at most $\maintenance - 1$ days before the end of the week.
The cost of $r$ is the sum of the costs of operating the legs $\ell$ in $r$ with airplane $j$.

The task consists in building, for each airplane $j$, a route satisfying the maintenance requirement in such a way that each leg of $\ell$ is operated by one airplane $j$, and the sum of the costs of the routes is minimum.

We now generalize our MIP to Air France tail assignment. 
Let $D' = (V,A')$ be the digraph with vertex set $V = \mathcal{L} \times \maintenance$, and arc set $A'$ composed of pairs $\big((\ell,\delta),(\ell',\delta')\big)$ such that $(\ell,\ell')$ is a tail assignment connection and $\ell$, $\delta$, $\ell'$, and $ \delta'$ satisfy one of the three conditions defining the arcs of digraph $D$ in Section~\ref{subsec:ip-ar}.
Digraph $D'$ is the analogue of digraph $D$ of Section~\ref{subsec:ip-ar} where airplane connections are replaced by tail assignment connections.
Contrary to digraph $D$, digraph $D'$ is acyclic as there is no connection between the end of the week and the beginning of the week.
For each airplane $j$, let $V^j$ and $ \tilde A^j$ be copies of $V$ and $A'$.
We build a digraph $D^j$ as follows.
Its vertex set is $V^j \cup \{s^j,t^j\}$, where $s^j$ is a source vertex, and $t^j$ a sink vertex.
Its arc set is denoted by $A^j$ and contains $\tilde A^j$ as well as arcs
\begin{itemize}
	\item $\big(s^j,(\ell,\delta)\big)$ such that leg $\ell$ starts from airport $k_0^j$ on day $\delta - \delta_0^j$ if $k_{0}^j$ is not a base,
	\item $\big(s^j,(\ell,1)\big)$ such that leg $\ell$ starts from airport $k_0^j$ if $k_{0}^j$ is a base,
	\item $\big((\ell,\delta),t^j\big)$ such that leg $\ell$ ends on day $d_{\ell}$ in an airport that is not a base, and $\delta < \maintenance + d_\ell - 8$,
	\item $\big((\ell,\delta),t^j\big)$ such that leg $\ell$ ends in a base.
\end{itemize}
The number $8$ in the third condition is just the number of days in a week plus $1$.
The digraph $D^j$ is acyclic.
For each leg $\ell$ and airplane $j$, we denote by $V_{\ell}^j$ the vertices of $V^j$ of the form $(\ell,\delta)$.
Given an arc $a$ in $A^j$, we define the cost $c_a$ to be equal to $0$ if the tail of $a$ is $s^j$ and to the cost of operating the leg of the tail vertex of $a$ with airplane $j$ otherwise. 
Given the definition of the digraph $D^j$, the following proposition is immediate.

\begin{prop}\label{prop:pathTA}
A sequence of legs $\ell_1,\ldots,\ell_k$ is a tail assignment route $r$ for airplane $j$ if and only if there exists $\delta_1,\ldots,\delta_k$ such that $s^j,(\ell_1,\delta_1),\ldots,(\ell_k,\delta_k),t^j$ is an $s^j$-$t^j$ path $P$ in $D^j$.
In that case, the cost of operating $r$  with $j$ is $\sum_{a \in P}c_a$.
\end{prop}

\noindent The following integer program therefore enables to model Air France tail assignment.

\begin{subequations}
        \makeatletter
        \def\@currentlabel{TA}
        \makeatother
\label{eq:ARprogramOpt}
        \renewcommand{\theequation}{TA.\arabic{equation}}
\begin{alignat}{2}
\min \enskip & \sum_{j \in [\na]} \sum_{a \in A^j} c_a x_a \\
\mathrm{s.t.}\enskip& \sum_{a\in\delta^+(s^j)}x_a = 1 && \qquad  \forall j \in [\na] \label{eq:aroptsource}\\
&\ds{\sum_{a\in\delta^-(v)}x_a = \sum_{a\in\delta^+(v)}x_a} && \qquad\forall v\in V^j, \forall j \in [\na] \smallskip \label{eq:aroptflow}\\
&\ds \sum_{j \in [\na]}{\sum_{a\in\delta^-(V_{\ell}^j)}x_a = 1} && \qquad\forall \ell\in\L \smallskip \label{eq:ARoptpartition} \\
&x_a \in \{0,1\} &&\qquad \forall a\in A^j, \forall j \in [\na]. \label{eq:ropt01}
\end{alignat}
\end{subequations}

The MIP proposed by Khaled et al.~for tail assignment with maintenance constraints is given by Equations (3-15) of their paper. 
They use binary variables $\overline{\bfx},\overline{\bfy}, \overline{\bfz}$. (We denote their variables with an overline to distinguish them from our variables.) The binary variables $\ovx_{ij}$ indicate if the leg $i$ is operated by the airplane $j$. The binary variables $\ovy_{jd}$ indicate if a maintenance of airplane $j$ takes place on day~$d$. The binary variables $\ovz_{ijd}$ indicate if a night maintenance of airplane $j$ takes place in the arrival airport of the leg $i$ on day $d$. 
Air France tail assignment version is slightly different from the one considered by Khaled et al. 
For instance,
they have a constraint limiting the cumulated flight time of an airplane between two maintenances (modeled by Equation~(13) of their MIP).
The following formulation adapts the MIP of Khaled et al. to Air France tail assignment problem.

\begin{subequations}
        \makeatletter
        \def\@currentlabel{KTA}
        \makeatother
\label{eq:TAkhaled}
        \renewcommand{\theequation}{KTA.\arabic{equation}}
\begin{alignat}{2}
\min_{(\ovbfx,\ovbfy,\ovbfz)}\enskip & \sum_{i \in F} \sum_{j \in P}c_{ij}\ovx_{ij} \\
\mathrm{s.t.}\enskip&\text{Equations~(3)-(6), (8), (9), (11), (12), (14), and (15) of \citep{khaled2018compact}}\\
& \sum_{\tilde d \in \{1,\ldots,\maintenance - \delta_0^j\}} y_{j\tilde d} \geq 1 \text{ for all $j$ in $[\na]$ } \label{eq:KTAbeginning}
\end{alignat}
\end{subequations}
where $F$ and $P$ are the notations of \citep{khaled2018compact} for the sets of flights legs and available airplanes respectively, and constraint~\eqref{eq:KTAbeginning} ensures that airplane $j$ spends its first night in a base after at most $\maintenance - \delta_0^j$ days.

The next proposition shows that the linear relaxation of Program~\eqref{eq:ARprogramOpt} provides bounds that will discard more partial solutions than that of ~\eqref{eq:TAkhaled} in a branch-and-bound.

\begin{prop}\label{prop:compareKhaled}
The optimal value of the linear relaxation of~\eqref{eq:ARprogramOpt} is not smaller than the one of the linear relaxation of~\eqref{eq:TAkhaled} and there are instances for which it is strictly larger.
\end{prop}

\begin{proof}

We first prove that any feasible solution of the linear relaxation of Program~\eqref{eq:ARprogramOpt} can be turned into a feasible solution of the linear relaxation of \eqref{eq:TAkhaled}. Consider a feasible solution $\bfx = (\bfx^1,\ldots,\bfx^{\na})$ of the linear relaxation of Program \eqref{eq:ARprogramOpt}, where $\bfx^j$ is the vector of variables $x_a$ with $a$ in $A^j$.
As Equation~\eqref{eq:aroptflow} defines the $s^j$-$t^j$ flow polyhedron of the acyclic digraph $D^j$, 
$\bfx^j$ can be written as a conic combination $\sum_{P}\lambda_P \chi^P$ of indicator vectors of $s^j$-$t^j$ paths in $D^j$, where the sum is taken over all $s^j$-$t^j$ paths and the $\lambda_P$ are non-negative.
Equation~\eqref{eq:aroptsource} gives $\sum_{P}\lambda_P = 1$.

Let $\calR^j$ be the set of tail assignment routes $r$ for airplane $j$. Given $r$ in $\calR^j$, we define $\lambda_r$ as the coefficient $\lambda_P$ of the $s^j$-$t^j$ path $P$ corresponding to $r$ according to Proposition~\ref{prop:pathTA}.
Denoting $c^r$ the cost of operating route $r$ with airplane $j$,
Proposition~\ref{prop:pathTA} also implies that
\begin{equation}\label{eq:lambdaCost}
 \sum_{a \in A^j}c_ax_a = \sum_{r \in \calR^j}\lambda_rc^r.
\end{equation}
Besides, as $\bfx$ satisfies Equation~\eqref{eq:ARoptpartition}, we have
\begin{equation}\label{eq:lambdaCover}
\sum_j \sum_{r \in \calR^{j,\ell}} \lambda_r = 1,
\end{equation}
 where $\calR^{j,\ell}$ denotes the subset of $\calR^j$ of routes containing leg $\ell$.

For each airplane $j$, route $r$ in $\calR^j$, leg $i$, and day $d$, we set $\ovx_{ij}^r = 1$ if leg $i$ is in route $r$ and $0$ otherwise, $\ovy_{jd}^r = 1$ if airplane $j$ operating route $r$ undergoes a maintenance on the night after day $d$ and $0$ otherwise, and $\ovz_{ijd}^r = 1$ if $i$ is in route $r$ and airplane $j$ operating $r$ undergoes a maintenance on the night after day $d$ and $0$ otherwise.
As $r$ satisfies the maintenance requirement, $(\ovx_{ij}^r,\ovy_{jd}^r,\ovz_{ijd}^r)_{id}$ satisfies the equations of the linear relaxation of \eqref{eq:TAkhaled} (restricted to airplane $j$) except cover constraint (3) of \citep{khaled2018compact}.
Furthermore, as $x_{ij}$ is the cost of operating leg $i$ with $j$, we have $\sum_{i \in F}c_{ij} \ovx_{ij}^r = c^r$, where $c^r$ is the cost of operating route $r$ with airplane $j$.

Let $(\ovbfx,\ovbfy,\ovbfz)$ be defined as follows.
For each airplane $j$, leg $i$, and day $d$, let  
$\ovx_{ij} = \sum_{r \in \calR^j} \lambda_r \ovx_{ij}^r$, let $\ovy_{jd} = \sum_{r \in \calR^j} \lambda_r \ovy_{jd}^r$, and let $\ovz_{ijd} = \sum_{r \in \calR^j} \lambda_r \ovz_{ijd}^r$. 
As $\lambda^r \geq 0$ and $\sum_{r \in \calR^j} \lambda_r = 1$, the components indexed by $j$ of $(\ovbfx,\ovbfy,\ovbfz)$ is a convex combination of the $(\ovx_{ij}^r,\ovy_{jd}^r,\ovz_{ijd}^r)_{id}$.
Hence, $(\ovbfx,\ovbfy,\ovbfz)$ satisfies the equations of the linear relaxation of \eqref{eq:TAkhaled} (restricted to airplane $j$) except cover constraint (3) of \citep{khaled2018compact}. 
Equation~\eqref{eq:lambdaCover} ensures that $(\ovbfx,\ovbfy,\ovbfz)$ also satisfies cover constraint (3) of \citep{khaled2018compact}, and is therefore a solution of the linear relaxation of \eqref{eq:TAkhaled}. Finally,
Equation~\eqref{eq:lambdaCost} ensures that $\sum_{i \in F} \sum_{j \in P}c_{ij}\ovx_{ij} = \sum_j \sum_{a \in A^j} c_a x_a$, which concludes the proof that the optimal value of the linear relaxation of~\eqref{eq:ARprogramOpt} is not smaller than the one of the linear relaxation of~\eqref{eq:TAkhaled}.

\begin{figure}
\begin{tikzpicture}
\def\l{0.5}
\def\h{1}
\node (a0) at (0.5*\l,0*\h) {};
\node (b0) at (2.5*\l,1*\h) {};
\node (c0) at (3*\l,1*\h) {};
\node (d0) at (5*\l,0*\h) {};
\node (e0) at (0*\l,0*\h) {};
\node (f0) at (5*\l,2*\h) {};
\node (g0) at (1.8*\l,2*\h) {};
\node (h0) at (6.8*\l,0*\h) {};

\draw[->,>=latex] (a0.center) to node[midway, right] {$\ell_a^1$} (b0.center);
\draw[->,>=latex] (c0.center) to node[midway, left ] {$\ell_b^1$} (d0.center);
\draw[->,>=latex] (e0.center) to[bend left] node[midway, left]  {$\ell_c^1$} (f0.center);

\node (a1) at (0.5*\l + 7*\l,0*\h) {};
\node (b1) at (2.5*\l + 7*\l,1*\h) {};
\node (c1) at (3*\l + 7*\l,1*\h) {};
\node (d1) at (5*\l + 7*\l,0*\h) {};
\node (e1) at (0*\l + 7*\l,0*\h) {};
\node (f1) at (5*\l + 7*\l,2*\h) {};
\node (g1) at (1.8*\l + 7*\l,2*\h) {};
\node (h1) at (6.8*\l + 7*\l,0*\h) {};

\draw[->,>=latex] (a1.center) to node[midway, right] {$\ell_a^2$} (b1.center);
\draw[->,>=latex] (c1.center) to node[midway, left ] {$\ell_b^2$} (d1.center);
\draw[->,>=latex] (e1.center) to[bend left] node[midway, left]  {$\ell_c^2$} (f1.center);
\draw[->,>=latex] (g1.center) to[bend left] node[midway, right] {$\ell_f^2$} (h1.center);

\node (g2) at (1.8*\l + 17*\l,2*\h) {};
\node (h2) at (6.8*\l + 17*\l,0*\h) {};

\draw[->,>=latex] (g2.center) to[bend left] node[midway, right] {$\ell_f^7$} (h2.center);


\draw[dashed] (-2*\l,0*\h) -- (14*\l,0*\h);
\draw[dashed] (-2*\l,1*\h) -- (14*\l,1*\h);
\draw[dashed] (-2*\l,2*\h) -- (14*\l,2*\h);

\draw[dashed] (15*\l,0*\h) -- (24*\l,0*\h) node[right] {Airport $A$ (base)};
\draw[dashed] (15*\l,1*\h) -- (24*\l,1*\h) node[right] {Airport $B$ (not-base)};
\draw[dashed] (15*\l,2*\h) -- (24*\l,2*\h) node[right] {Airport $C$ (not-base)};

\draw (-1*\l,-0.5*\h) -- (-1*\l,3*\h);
\draw (6*\l,-0.5*\h) -- (6*\l,3*\h);
\draw (13*\l,-0.5*\h) -- (13*\l,3*\h);
\draw (16*\l,-0.5*\h) -- (16*\l,3*\h);

\draw (23*\l,-0.5*\h) -- (23*\l,3*\h);

\node at (2*\l,2.5*\h) {Day 1};
\node at (9*\l,2.5*\h) {Day 2};
\node at (19*\l,2.5*\h) {Day 7};

\end{tikzpicture}
\bigskip

\begin{tabular}{|c|c|c|c|c|}
\hline
&\multicolumn{2}{c|}{Initial}&
\multicolumn{2}{c|}{Costs}\\
Airplane & airport $k_{0}^j$ & maint.~$\delta_0^j$ & $\ell_a$, $\ell_b$ & $\ell_c$,  $\ell_f$ \\ 
\hline
1 &A&1& 2 & 6 \\
2 &A&1& 3 & 9 \\
3 &A&1& 5 & 15 \\
\hline
\end{tabular}
\caption{Example used in the proof of Proposition~\ref{prop:compareKhaled}. }
\label{fig:khaled}
\end{figure}

Consider now the example with $n^{\mathrm{a}}=3$ and $\maintenance = 4$ on Figure~\ref{fig:khaled}. Dashed horizontal lines correspond to airports and arrows to legs between airports.
The weekly schedule is composed of six round trips between airport $A$ and airport $B$, and six round trips between airport $A$ and airport $C$.
There is an outward leg $\ell_a^d$ from $A$ to $B$ and an outward leg from $\ell_c^d$ from $A$ to $C$ every day $d$ in \{1,\ldots,6\}, a return leg $\ell_b^d$ from $B$ to $A$ every day $d$ in \{1,\ldots,6\}, and a return leg $\ell_f^d$ from $C$ to $A$ every day $d$ in \{2,\ldots,7\}.
The table provides, for each airplane $j$, the airport $k_{0}^j$ where the airplane starts, and the number of days $\delta_0^j$ since the last maintenance on day~1, and the costs of operating the different legs with the airplanes.
In the remaining of the proof, we show that the optimal value of the linear relaxation of~\eqref{eq:ARprogramOpt} on this instance is equal to 160, and exhibit a solution of the linear relaxation of $\eqref{eq:TAkhaled}$ with value 157.5.

Let $\bfx = (\bfx_j)_j$ be a solution of \eqref{eq:ARprogramOpt} on that instance.
Let $V_{\leq d}^j$ be the set vertices of $D^j$ composed of $s^j$ and vertices $(\ell_{\cdot}^{\tilde d},\delta)$ with $\tilde d \leq d$.
As the flow on the cut $\delta^+(V_{\leq d}^j)$ has value $1$, we have $\sum_{j=1}^3 \sum_{a \in \delta^+(V_{\leq d}^j)} x_a= 3$.
Besides, \eqref{eq:aroptflow} and \eqref{eq:ARoptpartition} give that 
$$\sum_{j=1}^3 \sum_{a \in \delta^+(V_{\ell_b^d}^j\cup V_{\ell_c^d}^j\cup V_{\ell_f^d}^j)}x_a = 3
\quad \text{and} \quad
\sum_{j=1}^3 \sum_{a \in \delta^-(V_{\ell_a^{d+1}}^j\cup V_{\ell_c^{d+1}}^j\cup V_{\ell_f^{d+1}}^j)}x_a = 3.
$$
Hence, any arc $a = (v,v')$ with $x_a > 0$ in $\delta^+(V_{\leq d}^j)$  is such that $v$ is in $ V_{\ell_b^d}^j\cup V_{\ell_c^d}^j\cup V_{\ell_f^d}^j$ and $v' $ belongs to $V_{\ell_a^{d+1}}^j\cup V_{\ell_c^{d+1}}^j\cup V_{\ell_f^{d+1}}^j$. 
Any arc $a$ with $x_a > 0$ is therefore of the form $\big((\ell,\delta)(\ell',\delta')\big)$ with $(\ell,\ell')$ in 
$$\big\{
(\ell_a^d,\ell_b^d),
(\ell_c^{d},\ell_f^{d+1}),
(\ell_b^{d},\ell_a^{d+1}),
(\ell_b^{d},\ell_c^{d+1}),
(\ell_f^{d},\ell_a^{d+1}),
(\ell_f^{d},\ell_c^{d+1})
\big\}.$$

Furthermore, as we discussed in the first part of the proof, $\bfx_j$ can be written as the conic combination $\sum_{r \in R^j} \lambda_r \chi^{P(r)}$, where $\chi^{P(r)}$ is the indicator vector of the $s^j$-$t^j$ path $P(r)$ corresponding to route $r$, and coefficients $\lambda_r$ are non-negative.
Remark that there is no connection $(\ell_f^d,\ell_c^{d+1})$ in a route that satisfies the maintenance requirement, as such a route would spend $\maintenance +1$ days out of a base.
Constraint \eqref{eq:ARoptpartition} applied to leg $\ell_c^{d+1}$ then ensures the only arcs $a$ such that $x_a > 0$ are of the form $\big((\ell,\delta),(\ell',\delta')\big)$ with $(\ell,\ell')$ in 
$$\big\{
(\ell_a^{d},\ell_b^{d}),
(\ell_b^{d},\ell_c^{d+1}),
(\ell_c^{d},\ell_f^{d+1}),
(\ell_f^{d},\ell_a^{d+1})
\big\}.$$
Hence, for each airplane $j$ the only routes $r$ in $\calR^j$  that can satisfy $\lambda_r > 0$ are 
$r_1 = \ell_a^1,
\ell_b^1,
\ell_c^2,
\ell_f^3,
\ell_a^4,
\ell_b^4,
\ell_c^5,
\ell_f^6$, 
$r_2 = \ell_c^1,
\ell_f^2,
\ell_a^3, 
\ell_b^3,
\ell_c^4,
\ell_f^5,
\ell_a^6,
\ell_b^6$, 
and 
$r_3 = \ell_a^2,
\ell_b^2,
\ell_c^3,
\ell_f^4,
\ell_a^5,
\ell_b^5,
\ell_c^6,
\ell_f^7$.
Operating any of these three routes has cost $32$ with airplane 1, $48$ with airplane $2$, and $80$ with airplane 3. 
Hence, $\bfx$ has cost $160$. 
Remark that this is the cost of the optimal integer solution obtained by assigning $r_j$ to airplane $j$.


On the contrary setting 

$$
\overline{x}_{\ell_a^d1} =
\overline{x}_{\ell_a^d2} =
\overline{x}_{\ell_b^d1} =
\overline{x}_{\ell_b^d2} =
\overline{z}_{\ell_b^d1d} =
\overline{z}_{\ell_b^d2d} =
\ovz_{\ell_b^61,7} =
\ovz_{\ell_b^62,7} =
\frac{1}{4} \text{ for all }d \in \{1,\ldots,6\}, 
$$
$$
\overline{x}_{\ell_c^d1} =
\overline{x}_{\ell_c^d2} = \frac{3}{8}
 \text{ for all }d \in \{1,\ldots,6\}, 
 \quad
\overline{x}_{\ell_f^d1} =
\overline{x}_{\ell_f^d2} =
\frac{3}{8} \text{ for all }d \in \{2,\ldots,7\} ,
$$
$$
\overline{x}_{\ell_a^d3} =
\overline{x}_{\ell_b^d3} =
\overline{z}_{\ell_b^d3d} =
\ovz_{\ell_b^63,7} =
\frac{1}{2} \text{ for all }d \in \{1,\ldots,6\}, 
$$
$$
\overline{x}_{\ell_c^d1} = \frac{1}{4} \text{ for all }d \in \{1,\ldots,6\},
\quad
\overline{x}_{\ell_f^d3} = 
\frac{1}{4} \text{ for all }d \in \{2,\ldots,7\},
$$
$$
\overline{y}_{1d}       =
\overline{y}_{2d}       =
\frac{1}{4} \text{ for all }d \in \{1,\ldots,7\},
\quad \text{and}\quad
\overline{y}_{3d} = \frac{1}{2} \text{ for $d$ in $\{1,\ldots,7\}$},
$$
and $z_{ijd} = 0$ for any $i,j,d$ such that $z_{ijd}$ has still not been defined, provides a feasible solution of the linear relaxation of \eqref{eq:TAkhaled} with cost 157.5, which concludes the proof.
\end{proof}

\section{Example of the Monoid resource constrained shortest path problem}
\label{sec:example_of_the_monoid_resource_constrained_shortest_path_problem}
\label{sec:appendix_example}

\begin{figure}
\begin{tabular}{lc}

a. &
\begin{tikzpicture}
\def\l{0.8}
\def\h{1}

\draw (3*\l,-0.5*\h) -- (3*\l,2.5*\h);
\draw (10*\l,-0.5*\h) -- (10*\l,2.5*\h);

\node at (1*\l,2.5*\h) {Day $i-1$};
\node at (6*\l,2.5*\h) {Day $i$};
\node at (12*\l,2.5*\h) {Day $i+1$};

\draw[->,>=latex] (0*\l,-1*\h) -- (18*\l,-1*\h) node[right,below] {time};

\draw[dashed] (0*\l,2*\h) -- (14*\l,2*\h) node[right] (a1) {Airport $A_1$};
\draw[dashed] (0*\l,1*\h) -- (14*\l,1*\h) node[right] (a2) {Airport $A_2$};
\draw[dashed] (0*\l,0*\h) -- (14*\l,0*\h) node[right] (a3) {Airport $A_3$};

\node[right = 0 of a1] (b1) {};
\node[right = 0 of a2] (b2) {};
\node[right = 0 of a3] (b3) {};

\draw[<->] (b1) to node[midway,right] {6h flight} (b2);
\draw[<->] (b2) to node[midway,right] {2h flight} (b3);

\node (l1t) at (0*\l,2*\h) {};
\node (l1h) at (2.3*\l,1*\h) {};
\node (l2t) at (1*\l,0*\h) {};
\node (l2h) at (2*\l,1*\h) {};
\node (l3t) at (4*\l,1*\h) {};
\node (l3h) at (5*\l,0*\h) {};
\node (l4t) at (5*\l,0*\h) {};
\node (l4h) at (6*\l,1*\h) {};
\node (l5t) at (7*\l,1*\h) {};
\node (l5h) at (9.3*\l,2*\h) {};
\node (l6t) at (7*\l,1*\h) {};
\node (l6h) at (8*\l,0*\h) {};
\node (l7t) at (8*\l,0*\h) {};
\node (l7h) at (9*\l,1*\h) {};
\node (l9t) at (10.7*\l,2*\h) {};
\node (l9h) at (13*\l,1*\h) {};
\draw[->,>=latex] (l1t.east) to node[color=black, midway, right] {$\ell_{1}$} (l1h.west);
\draw[->,>=latex] (l2t.east) to node[midway, right] {$\ell_{2}$} (l2h.west);
\draw[->,>=latex] (l3t.east) to node[color=black, midway, right] {$\ell_{3}$} (l3h.west);
\draw[->,>=latex] (l4t.east) to node[color=black, midway, right] {$\ell_{4}$} (l4h.west);
\draw[->,>=latex] (l5t.east) to node[midway, right] {$\ell_{5}$} (l5h.west);
\draw[->,>=latex] (l6t.east) to node[midway, right] {$\ell_{6}$} (l6h.west);
\draw[->,>=latex] (l7t.east) to node[midway, right] {$\ell_{7}$} (l7h.west);
\draw[->,>=latex] (l9t.east) to node[midway, right] {$\ell_{8}$} (l9h.west);

\end{tikzpicture}
\bigskip
\bigskip
 \\
 b.&

\begin{tikzpicture}
\def\l{2}
\def\h{1.5}

\node[draw,circle] (v1) at (-1*\l,2*\h) {$\ell_1$};
\node[draw,circle] (v2) at (-1*\l,0*\h) {$\ell_2$};
\node[draw,circle] (v3) at (0.7*\l,0*\h) {$\ell_3$};
\node[draw,circle] (v4) at (2.1*\l,0*\h) {$\ell_4$};
\node[draw,circle] (v5) at (3*\l,2*\h) {$\ell_5$};
\node[draw,circle] (v6) at (3.5*\l,0*\h) {$\ell_6$};
\node[draw,circle] (v7) at (4.9*\l,0*\h) {$\ell_7$};
\node[draw,circle] (v9) at (5.2*\l,2*\h) {$\ell_8$};

\draw[<-,>=latex,dashed] (v1) -- ++(-1,0);
\draw[<-,>=latex,dashed] (v2) -- ++(-1,0);
\draw[->,>=latex,dashed] (v7) -- ++(1,0);
\draw[->,>=latex,dashed] (v9) -- ++(1,0);

\draw[arc] (v1) to[color=black] node[midway,left] {$\big((0,0,2,2),-6.2\big)$} (v3);
\draw[arc] (v1) to node[midway,above] {$\big((0,0,1,6),6.2\big)$} (v5);
\draw[arc] (v1) to node[midway,near start] {$\big((0,0,1,2),2.1\big)$} (v6);
\draw[arc] (v2) to  node[midway,below] {$\big((0,0,1,2),4.9\big)$} (v3);
\draw[arc] (v2) to  node[near end] {$\big((0,0,1,6),6.2\big)$} (v5);
\draw[arc] (v2) to[bend right] node[midway,below] {$\big((0,0,1,2),5.4\big)$} (v6);
\draw[arc] (v3) to[color=black]  node[midway,below] {$\big((1,2),-2.1\big)$} (v4);
\draw[arc] (v4) to node[midway,right] {$\big((1,6),1.3\big)$} (v5);
\draw[arc] (v4) to node[midway,below] {$\big((1,2),2.1\big)$} (v6);
\draw[arc] (v6) to node[midway,below] {$\big((1,2),1.5\big)$} (v7);
\draw[arc] (v5) to node[midway,above] {$\big((0,0,2,6),1.2\big)$} (v9);
\end{tikzpicture}
\end{tabular}
\caption{Example of instance of the monoid shortest path problem}
\label{fig:exampleMonoidGraph}
\end{figure}

This appendix details the execution of Algorithm~\ref{alg:enumeration} on a simple example.
Figure~\ref{fig:exampleMonoidGraph} provides an example of instance of the crew pairing pricing subproblem and its \MRCSP modeling. Only a subpart of the instance is represented. 
On Figure~\ref{fig:exampleMonoidGraph}.a, legs are represented as arrows between airports. 
On this instance, there are three airports $A_1$, $A_2$, and $A_3$, and legs only between $A_1$ and $A_2$, and $A_2$ and $A_3$. The flying durations between $A_1$ and $A_2$ and $A_2$ and $A_3$ are respectively 6 and 2 hours. There are two reduced rests: $(\ell_1,\ell_3)$ and $(\ell_5,\ell_8)$. The maximum duty flying duration $F(t)$ in a duty is taken equal to $F_m=9$ hours for all $t$.

Figure~\ref{fig:exampleMonoidGraph}.b provides the corresponding digraph $D$, as well as the resource of each arc.
Note that the component in $M^\rmeas$ of the resources of day connections is in $\Z_+ \times \R_+$ and that of the resources of night connections is in $(\Z_+ \times \R_+)^2$. Furthermore, reduced rests are the only night connections whose resources have a component in $M^\rmeas$ of the form $(0,0,2,\cdot)$, the other night connections having a component of the form $(0,0,1,\cdot)$. 
The reduced costs have been chosen arbitrarily.
For simplicity, we assume that $\big((0,0,0,0),0\big)$ is a lower bound on the resource of $\ell_7$-$d$ paths and on the resource of $\ell_8$-$d$ paths.

\subsubsection*{Algorithm~\ref{alg:enumeration} execution} 
\label{ssub:algorithm_execution}

\begin{table}
\begin{tabular}{|c|c|}
\hline
$v$ 
& $B_v$ \\
\hline
$\ell_1$ & \big\{\big((0,0,0,0),-5.8\big)\big\}\\
$\ell_3$ & \big\{\big((2,6,0,0),0.4\big)\big\}\\
$\ell_4$ & \big\{\big((1,4,0,0),2.5\big)\big\}\\
$\ell_5$ & \big\{\big((0,0,0,0),1.2\big)\big\}\\
$\ell_6$ & \big\{\big((1,2,0,0),1.5\big)\big\}\\
\hline
\end{tabular}
\caption{Sets of bounds $B_v$ used (here, singletons)}
\label{tab:bounds}
\end{table}

We now provide an example of Algorithm~\ref{alg:enumeration} when the bounds in Table~\ref{tab:bounds} are used.
We remind the reader that these bounds are computed before the execution of Algorithm~\ref{alg:enumeration} in a preprocessing and taken by this latter in input. 
See Section~\ref{sub:enhanced_bounds_computations_for_the_column_generation_context} for more details on bounds. 
At the very end of the appendix, we illustrate the way bounds are computed by justifying the set of bounds $B_{\ell_3}$ in this table.

Let $P_1$ be an $o$-$\ell_1$ path with resource $\big((0,0,1,6),1.0\big)$.
Table~\ref{tab:singleBounds} describes the iterations of Algorithm~\ref{alg:enumeration} where $P_1$ and the paths starting by $P_1$ are dealt with.
Each iteration is separated by an horizontal line.
Column $P$ provides the path $P$ considered at Step~\ref{step:pathSelection} of Algorithm~\ref{alg:enumeration}.
We assume that $\mathsf{L}_{\ell_3}^{\mathrm{nd}}$, $\mathsf{L}_{\ell_5}^{\mathrm{nd}}$, and $\mathsf{L}_{\ell_6}^{\mathrm{nd}}$ are empty when $P_1$ is considered as path $P$, and  $\mathsf{L}_{\ell_4}^{\mathrm{nd}}$ is empty when $P_1,\ell_3$ is considered.
We also assume that $c_{od}^{UB} = +\infty$ and hence $c(q_P\rplus b) \leq c_{od}^{UB}$ during all the iterations detailed.
Column key$(P)$ provides its key defined in Equation~\eqref{eq:keyDefinition}. 
As the treatment of reduced costs is standard,  to enhance readability, we omit them in all resources in Table~\ref{tab:singleBounds} and in the remaining of the discussion.
Column $q_P$ gives the resource of $P$, column $a$ provides the arc of Step~\ref{step:arcLoop}, column $q_a$ gives its resource.
Path $Q$ of Step~\ref{step:Qdef} is path $P$ followed by $a$. Then next column gives the resource $q_Q = q_P \rplus q_a$ of $Q$, and column $w$ provides the destination of $Q$ computed at Step~\ref{step:Qstarts}. Column $b$ provides the single bound in $B_w$, and the next column provides $q_Q \rplus b$ computed at Step~\ref{step:LBtest}. 
Finally, the last column indicates if $Q$ is added to $\mathsf{L}$ at Step~\ref{step:updateL}.

The key of $P_1$ is equal to $-4.8$ because, with $P= P_1$ and $B_{\ell_1}= \{b\}$, we have $q_P \rplus b = \big((0,0,1,6),1.0\big) \rplus \big((0,0,0,0),-5.8\big) ) = \big((0,0,0,0),-4.8\big)$, and $c\big((0,0,0,0),-4.8\big) = -4.8$.
The sums $(0,0,2,2)\rplus(2,6,0,0)$ and $(0,0,3,4) \rplus (1,4,0,0)$ are equal to $(0,0,0,0)$ because $(2,2) + (2,6) = (3,4) +(1,4) = (4,8) \leq (4,F_m)$, where $F_m = 9$. 
Since $\rho((0,0,0,0))=0$, the path $Q$ cannot be discarded at the iterations where $P = P_1$ and $a=(\ell_1,\ell_3)$,  and where $P=P_1,\ell_3$ and $a=(\ell_3,\ell_4)$, and it is added to $\mathsf{L}$.  
On the contrary $(0,0,4,10) \rplus (0,0,0,0) = \infty$ because $(4,10) \nleq (4,F_m)$. 
Since $\rho(\infty)=1$, the path $Q$ is not kept after Step~\ref{step:LBtest} when $a=(\ell_4,\ell_5)$, and it is not added to $\mathsf{L}$.
We have a similar outcome when $a=(\ell_4,\ell_6)$: In this case, $(0,0,4,6) \rplus (1,2,0,0) = \infty$ because $(5,8) \nleq (4,F_m)$.
The treatment of $Q= P_1,\ell_3,\ell_4,\ell_6$ shows the interest of the bounds: Although path $Q$ itself satisfies rule~\ref{rule:NbLegs},
the algorithm identifies that any path starting by $Q$ violates rule~\ref{rule:NbLegs}.

\newcommand{\specialcell}[2][c]{%
  \begin{tabular}[#1]{@{}c@{}}#2\end{tabular}}

\begin{table}
\begin{outdent}

\begin{tabular}{|c|c|c|c|c|c|c|c|c|c|}
\hline
\specialcell{$P$} &
\specialcell{key$(P)$} &
\specialcell{$q_P$}&
\specialcell{$a$} &
\specialcell{$q_a$}&
\specialcell{$q_Q$} &
\specialcell{$w$ } &
\specialcell{$b$ \\ $B_w= \{b\}$} &
\specialcell{$q_Q \rplus b$ \\} &
\specialcell{$Q$ added \\ to $\mathsf{L}$} \\
\hline
$P_1$ &
$-4.8$ &
$(0,0,1,6)$&
$(\ell_1,\ell_3$) &
$(0,0,2,2)$&
$(0,0,2,2)$ &
$\ell_3$ &
$(2,6,0,0)$ &
$(0,0,0,0)$ &
\textsf{yes}
\\
&&&
$(\ell_1,\ell_5$) &
$(0,0,1,6)$&
$(0,0,1,6)$ &
$\ell_5$ &
$(0,0,0,0)$ &
$(0,0,0,0)$ &
\textsf{yes}
\\
&&&
$(\ell_1,\ell_6$) &
$(0,0,1,2)$&
$(0,0,1,2)$ &
$\ell_6$ &
$(1,2,0,0)$ &
$(0,0,0,0)$ &
\textsf{yes}
\\\hline
$P_1,\ell_3$ &
$-4.8$ &
$(0,0,2,2)$&
$(\ell_3,\ell_4)$&
$(1,2)$&
$(0,0,3,4)$ & 
$\ell_4$ &
$(1,4,0,0)$&
$(0,0,0,0)$&
\textsf{yes}
\\\hline
$P_1,\ell_3,\ell_4$ &
$-4.8$ &
$(0,0,3,4)$ & 
$(\ell_4,\ell_5)$&
$(1,6)$&
$(0,0,4,10)$&
$\ell_5$&
$(0,0,0,0)$&
$\infty$&
\textsf{no}
\\
&&&
$(\ell_4,\ell_6)$&
$(1,2)$&
$(0,0,4,6)$&
$\ell_6$&
$(1,2,0,0)$&
$\infty$&
\textsf{no}
\\\hline
$P_1,\ell_6$ &
$4.6$&
$(0,0,1,2)$&
$(\ell_6,\ell_7)$&
$(1,2)$&
$(0,0,2,4)$ & 
$\ell_7$ &
\multicolumn{3}{c|}{\multirow{2}{*}{etc.}}
\\\cline{1-7}
$P_1,\ell_5$ &
$8.4$&
$(0,0,1,6)$&
$(\ell_5,\ell_8)$&
$(0,0,2,6)$&
$(0,0,2,6)$ & 
$\ell_8$ &
\multicolumn{3}{c|}{}
\\
\hline
\end{tabular}
\end{outdent}
\caption{Algorithm execution with single bounds: iterations considering paths $P$ starting by an $o$-$\ell_1$ path $P_1$}
\label{tab:singleBounds}
\end{table}

\begin{table}
\begin{outdent}

\begin{tabular}{|c|c|c|c|c|c|c|c|c|c|}
\hline
\specialcell{$P$} &
\specialcell{key$(P)$} &
\specialcell{$q_P$}&
\specialcell{$a$} &
\specialcell{$q_a$}&
\specialcell{$q_Q$} &
\specialcell{$w$ } &
\specialcell{$b$} &
\specialcell{$q_Q \rplus b$ \\} &
\specialcell{$Q$ added \\ to $\mathsf{L}$} \\
\hline
\specialcell{$P_1$ }&
$-4.8$&
{$(0,0,1,6)$}&
{$(\ell_1,\ell_3$) }&
{$(0,0,2,2)$}&
{$(0,0,2,2)$ }&
{$\ell_3$ }&
$(2,8,0,0)$   &
$\infty$  &
\\ 
&&&&&&&
$(3,6,0,0)$& $\infty$& \textsf{no} 
\\
&&&
$(\ell_1,\ell_5$) &
$(0,0,1,6)$&
$(0,0,1,6)$ &
$\ell_5$ &
$(0,0,0,0)$ &
$(0,0,0,0)$ &
\textsf{yes}
\\
&&&
$(\ell_1,\ell_6$) &
$(0,0,1,2)$&
$(0,0,1,2)$ &
$\ell_6$ &
$(1,2,0,0)$ &
$(0,0,0,0)$ &
\textsf{yes}
\\\hline
$P_1,\ell_6$ &
$4.6$&
$(0,0,1,2)$&
$(\ell_6,\ell_7)$&
$(1,2)$&
$(0,0,2,4)$ & 
$\ell_7$ &
\multicolumn{3}{c|}{\multirow{2}{*}{etc.}}
\\\cline{1-7}
$P_1,\ell_5$ &
$8.4$&
$(0,0,1,6)$&
$(\ell_5,\ell_8)$&
$(0,0,2,6)$&
$(0,0,2,6)$ & 
$\ell_8$ &
\multicolumn{3}{c|}{}
\\
\hline
\end{tabular}
\end{outdent}
\caption{Algorithm execution with $B_{\ell_3} = \big\{\big((2,8,0,0),0.4\big),\big((3,6,0,0),1.5\big)\big\}$}
\label{tab:multiBounds}
\end{table}


Table~\ref{tab:multiBounds} provides the same informations as Table~\ref{tab:singleBounds} when we use a set of bounds  $B_{\ell_3} = \big\{\big((2,8,0,0),0.4\big),\big((3,6,0,0),1.5\big)\big\}$ instead of the singleton given in Table~\ref{tab:bounds}. 
Fewer iterations are then needed: 
the ``if'' condition at Step~\ref{step:LBtest} is not satisfied when $P_1,\ell_3$ is considered as path $Q$, and path $P_1,\ell_3$ is never added to $\mathsf{L}$. 
Even though $P_1,\ell_3$ itself satisfies rules~\ref{rule:NbLegs} and~\ref{rule:FlyingTime}, the algorithm identifies that any $o$-$d$ path starting by $P_1,\ell_3$ does not satisfy at least one of these rules.

\subsubsection*{Rationale of $B_{\ell_3}$}
We explain why $B_{\ell_3}$ is a correct bound set, both in the singleton and non-singleton cases. This explanation can be seen as a rough illustration of the procedure mentioned in Section~\ref{sub:enhanced_bounds_computations_for_the_column_generation_context} for the bound computation.

Any $\ell_3$-$d$ path must either start with $\ell_3,\ell_4,\ell_5,\ell_8$, or with $\ell_3,\ell_4,\ell_6,\ell_7$. 
Recall that we have assumed that $\big((0,0,0,0),0\big)$ is a lower bound on the resource of $\ell_7$-$d$ paths and on the resource of $\ell_8$-$d$ paths.
Given that $\re_{(\ell_3,\ell_4)} \rplus \re_{(\ell_4,\ell_5)} \rplus \re_{(\ell_4,\ell_8)} \rplus \big((0,0,0,0),0\big) = \big((2,8,0,0),0.4\big)$ and 
$\re_{(\ell_3,\ell_4)} \rplus \re_{(\ell_4,\ell_6)} \rplus \re_{(\ell_6,\ell_7)} \rplus \big((0,0,0,0),0\big) = \big((3,6,0,0),1.5\big)$, any $\ell_3$-$d$ path starting by $\ell_3,\ell_4,\ell_5$ has a resource lower bounded by $\big((2,8,0,0),0.4\big)$, and any $\ell_3$-$d$ path starting by $\ell_3,\ell_4,\ell_6,\ell_7$ has a resource lower bounded by $\big((3,6,0,0),1.5\big)$. This explains why $\big\{\big((2,8,0,0),0.4\big),\big((3,6,0,0),1.5\big)\big\}$ can be used as a set of bounds $B_{\ell_3}$, and why $\big((2,8,0,0),0.4\big)\meet\big((3,6,0,0),1.5\big) =\big((2,6,0,0),0.4\big)$ is a lower bound on the resource of any $\ell_3$-$d$ path.


\end{document}